\documentclass[a4paper]{article}

\usepackage{amsfonts}
\usepackage{latexsym}
\usepackage{amssymb}
\usepackage{amscd}

 \usepackage{amsfonts}
\usepackage{latexsym}
\usepackage{amssymb}
\usepackage{amscd}
  
 \newtheorem{thm}{Theorem}[section]

\newtheorem{defn}[thm]{Definition}

\newtheorem{lemma}[thm]{Lemma}
\newtheorem{prop}[thm]{Proposition}
\newtheorem{cor}[thm]{Corollary}

\newenvironment{prf}{{\bf Proof.}}{\hfill$\Box$\\[1mm]}
 
 \newcommand{\A}{{\mathcal{A}}}

\newcommand{\I}{{\mathcal{I}}}

\def\a{\alpha} 

\def\b{\beta}

\def\d{\delta}
\def\D{\Delta}
\def\e{\varepsilon}

\def\l{\lambda}

\def\t{\tau}
\def\om{\omega} 
\def\Om{\Omega}

\def\aap{a'}
\def\ep{e'}

\def\bC{\mathbb{C}}
\def\bN{\mathbb{N}} 
\def\bQ{\mathbb{Q}}
\def\bR{\mathbb{R}}  
 
\def\bZ{\mathbb{Z}}

\def\A{\mathcal{A}}
\def\C{\mathcal{C}}
 
\def\F{\mathcal{F}}
\def\H{\mathcal{H}}

\def\K{\mathcal{K}}

\def\M{\mathcal{M}}

 \def\S{\mathcal{S}}

\def\id{\mathrm{id}}
 
\def\tr{\mathrm{tr}}

\def\assets{\mathrm{as\;sets}}

\def\other{\mathrm{other}}

\def\otherwise{\mathrm{otherwise}}

\def\field{k}
\def\Heis{\mathrm{Heis}}

\def\half{{\frac{1}{2}}}

\def\br{{\bf r}}
\def\tt{\tilde{t}}

\def\hhcm{U({\bf b_+})}
\def\hhcmdual{\field [ \Bplus ]}

\def\ahcm{\field [ \Dplus ] }
\def\ahcmdual{U(\dplus)}

\def\hcm{\ahcm \lrbicross \hhcm}
\def\hcmdual{\ahcmdual \rlbicross\hhcmdual}
\def\ulh{U_\lambda({\bf heis})}
\def\extendedulh{{\tilde{U_\lambda}}({\bf heis})}

\def\clh{\field_\l [ \Heis]}
\def\ulbplus{U_\l({\bf b_+})}

\def\a{\alpha}
\def\b{\beta}

\def\rto{\rightarrow}

\def\aap{a'}
\def\cp{c'}
\def\ep{e'}
\def\fp{f'}
\def\gp{g'}
\def\hp{h'}
\def\kp{k'}

\def\tr{\triangleright}
\def\tl{\triangleleft}

\def\HCM{\H_{\mathrm{CM}}}

\def\HCMoriginal{\H_{\mathrm{CM}}^{\mathrm{left}}}
\def\HCMlambda{\H_{\mathrm{CM}}^\lambda}
\def\UCM{U_{\mathrm{CM}}}
\def\UCMlambda{U_{\mathrm{CM}}^\lambda}

\def\Diff{\mathrm{Diff}}

\def\nchoosej{\left(\begin{array}{cr} n\\ j \end{array} \right)}

\def\nplusonechoosek{\left(\begin{array}{cr} n +1\\ k \end{array} \right)}

\def\nminusonechoosekminusone{\left(\begin{array}{cr} n-1\\ k-1 \end{array} \right)}
\def\nminusonechoosek{\left(\begin{array}{cr} n-1\\ k \end{array} \right)}
\def\rchooses{\left(\begin{array}{cr} r\\ s \end{array} \right)}

\def\bplus{\mathrm{{\bf b}_+}}
\def\Bplus{\mathrm{B_+}}
\def\Dplus{{\mathrm{D}_0}}
\def\dplus{{\bf d}_0}
\def\FBplus{F[ \Bplus]}
\def\Ulbplus{U_\l ({\bf b_+})}
\def\kt{\field[t]}
\def\hbg{{\mathbb{H}}_3 (\field)}

\def\bg{\bf g}

\def\Adual{\A'}
\def\Hdual{\H'}

\def\coleft{{(0)}}
\def\comid{{\overline{(1)}}}
\def\coright{{(2)}}

\def\oneb{{\overline{(1)}}}
\def\twob{{\overline{(2)}}}
\def\one{{(1)}}
\def\two{{(2)}}

\def\modulostuff{\mathrm{modulo}\; \HCMlambda \otimes \I + \I \otimes \HCMlambda}

\def \sltwoR{{\bf sl_2}(\bR)}
\def\SLtwoR{SL_2 (\bR)}

\newcommand{\la}{\triangleright}
\newcommand{\ra}{\triangleleft}
\newcommand{\del}{{\partial}}

\newcommand{\lcross}{{>\!\!\!\triangleleft}}

\newcommand{\rlbicross}{{\triangleright\!\!\!\blacktriangleleft}}
\newcommand{\lrbicross}{{\blacktriangleright\!\!\!\triangleleft}}

\newcommand{\tens}{{\otimes}}

\begin{document}

\title{Bicrossproduct approach to the Connes-Moscovici Hopf algebra}
\author{Tom~Hadfield\footnote{Supported by an EPSRC postdoctoral fellowship} ,~Shahn~Majid\footnote{Supported during completion of the work by the Perimeter Institute, Waterloo, Ontario.}}
\date{\today}
\maketitle

\centerline{School of Mathematical Sciences,}
\centerline{Queen Mary, University of London}
\centerline{327 Mile End Road, London E1 4NS, England}
\centerline{t.hadfield@qmul.ac.uk, s.majid@qmul.ac.uk}
\centerline{}
\centerline{MSC (2000): 58B32, 58B34, 17B37}

   \abstract{We give a rigorous proof that the (codimension one) Connes-Moscovici Hopf algebra $\HCM$ is isomorphic to a bicrossproduct Hopf algebra linked to a group 
    factorisation of the diffeomorphism group $\Diff^+ (\bR)$.
     We construct a second bicrossproduct $\UCM$ equipped with a nondegenerate dual pairing with $\HCM$. 
      We give a  natural  quotient Hopf algebra $\clh$ of $\HCM$ and Hopf subalgebra $\ulh$ of $\UCM$ which again are in duality. All these Hopf algebras arise as deformations of commutative or cocommutative Hopf algebras that we describe in each case. Finally we develop the noncommutative differential geometry of $\clh$ by studying first order differential calculi of small dimension.
 } 

\section{The Connes-Moscovici Hopf algebra $\HCM$}
\label{introduction}

The Connes-Moscovici Hopf algebras originally appeared in \cite{cm}, arising  from a longstanding internal problem of noncommutative geometry,  the computation of the index of transversally elliptic operators on foliations. This family of Hopf algebras (one for each positive integer) was found to reduce transverse geometry to a universal geometry of affine nature, and provided the initial impetus for the development of Hopf-cyclic cohomology. 
The cyclic cohomology of these Hopf algebras was shown by Connes and Moscovici to serve as an organizing principle for the computation of the cocycles in their local index formula \cite{cm95}.
They are also closely related to the Connes-Kreimer Hopf algebras of rooted trees arising from renormalization of quantum field theories \cite{ck1}. More recently  these Hopf algebras have appeared in number theory,  in the context of operations on spaces of modular forms and modular Hecke algebras \cite{cm04a} and spaces of $\bQ$-lattices \cite{cmar}. They appear to play a near-ubiquitous role as symmetries in noncommutative geometry. 
There is also an algebraic approach to diffeomorphism groups \cite{SM_qbdg}, which we link to Connes and Moscovici's work.

In this paper we focus on the simplest example, the codimension one Connes-Moscovici Hopf algebra.
 We work with a right-handed version of this algebra, which we denote $\HCM$.  The algebras in  \cite{cm} were implicitly defined over $\bR$ or $\bC$,  but throughout this paper we will work over an arbitrary field $\field$ of characteristic zero.

\begin{defn} We define $\HCM$ to be the Hopf algebra (over $\field$)  generated by elements $X$, $Y$, $\d_n$ 
($n \geq 1$), with 
\begin{eqnarray}
\label{HCM_coproduct}
&&[Y,X] = X, \quad [X, \d_n ] = \d_{n+1} , \quad [Y, \d_n ] = n \d_n, \quad [ \d_m , \d_n ] =0 \quad  \forall \; m,n\nonumber\\
&&\D(X)=X \otimes 1 + 1 \otimes X + Y \otimes \d_1,\nonumber\\   
&&\D(Y)=Y \otimes 1 + 1 \otimes Y,\quad \D( \d_1 )= \d_1 \otimes 1 + 1 \otimes \d_1\nonumber\\
&&\e (X) = 0 = \e(Y), \quad  \e( \d_n)=0 \quad \forall \; n\nonumber\\
&&S(Y) = -Y,  \quad S(X) = Y \d_1 -X, \quad S( \d_1) = - \d_1
\end{eqnarray}
with $\D ( \d_{n+1} )$, $S(\d_{n+1})$ defined inductively from the relation $[X, \d_n ] = \d_{n+1}$.
\end{defn}

This differs from the Hopf algebra defined in \cite{cm}, p206, in that Connes and Moscovici take $\D(X) = X \otimes 1 + 1 \otimes X + \d_1 \otimes Y$. 
 We will denote this original left-handed version by $\HCMoriginal$.

The first part of this paper gives a rigorous proof that $\HCM$ is isomorphic to  a bicrossproduct Hopf algebra linked to a factorisation of the group $\Diff^+ (\bR)$ 
of positively-oriented diffeomorphisms of the real line. 
Recall that a group $X$ is said to factorise into subgroups $G$ and $M$  if group multiplication gives a set bijection $G \times M \rightarrow X$. 
 We write $X = G \bowtie M$. 
  As remarked in \cite{cm},  the group
  \begin{equation}\label{defn_Diffplus}
  \Diff^+ (\bR) = \{ \; \varphi \in \Diff (\bR) \; : \; \varphi'(x) > 0 \; \; \forall \; x \in \bR \; \}
  \end{equation}
 factorises into the two subgroups
\begin{equation}\label{defn_Dzero}
\Dplus = \Diff^+_0 (\bR) = \{ \; \phi \in \Diff^+ (\bR) \; : \; \phi(0) =0, \; \phi'(0) =1 \; \}
\end{equation}
 and 
$\Bplus = \{ \; (a,b) : x \mapsto ax+b \; : \; a, \, b \in \bR, \; a>0 \; \}$
the subgroup of affine diffeomorphisms, which 
 we identify with its faithful matrix  representation
\begin{equation}\label{defn_Bplus}
\Bplus = \{ (a,b) =  {\tiny\left(
\begin{array}{cc}
a & b \\
0 & 1
\end{array}
 \right)} : \; a,b \in \bR, \; a > 0 \; \}
\end{equation}

Given a group factorisation $X = G \bowtie M$ of a finite group $X$,  there is a natural construction of dually-paired 
 finite-dimensional bicrossproduct Hopf algebras denoted $k[M] \lrbicross kG$,  $kM \rlbicross k[G]$ \cite{kac,SM:thesis,takeuchi}.
In our case, the group $X = \Diff^+ (\bR)$ is very far from finite, so it remains a challenge to construct an analogous pair of infinite-dimensional  bicrossproduct Hopf algebras.
 Fortunately the  bicrossproduct construction is more general than the group factorisation case.
Given Hopf algebras $\A$ and $\H$, with $\A$ a left $\H$-module algebra and $\H$ a right $\A$-comodule coalgebra with action and coaction compatible in an appropriate sense, then it is possible to equip the vector space $\A \otimes \H$ with the structure of a Hopf algebra, the (left-right) bicrossproduct denoted $\A \lrbicross \H$ \cite{SM_hvn,SM_book}.
 Similarly, we can construct a (right-left) bicrossproduct 
 $\H \rlbicross \A$ from a right $\H$-module algebra $\A$ and  left $\A$-comodule coalgebra $\H$. 
 Thus in our  case we define Hopf algebras $\ahcm$, $\hhcm$, together with a  left action of $\hhcm$ on $\ahcm$ and  right coaction of $\ahcm$ on $\hhcm$ which  we prove  are compatible in the sense necessary for the construction of a bicrossproduct Hopf algebra $\hcm$ (Theorem \ref{compatible}).  
 We  prove that $\hcm$ is  isomorphic to  $\HCM$.
 We then construct a second bicrossproduct Hopf algebra $\UCM = \hcmdual$ (Proposition \ref{dual_of_HCM})  equipped with a nondegenerate dual pairing with $\HCM$.

We explain carefully how the actions and coactions giving rise to these bicrossproducts can be derived from the factorisation $\Diff^+ (\bR) = \Bplus \bowtie \Dplus$.  
 This  serves as motivation and is not part of our proof. 
 However, if we  simply presented compatible actions and coactions without indicating how they arose, although no rigour would be lost this would leave things very opaque. 
 
 We note that the original Connes-Moscovici Hopf algebra $\HCMoriginal$, defined as in (\ref{HCM_coproduct}), but with $\d_1 \otimes Y$ rather than $Y \otimes \d_1$ appearing in $\D(X)$, can also be shown to be a bicrossproduct linked to this group factorisation.
 The construction is given in Section \ref{section_HCM_as_bicross}. 
 However, as we explain $\HCM$ rather than $\HCMoriginal$ is in an appropriate sense the natural bicrossproduct associated to this factorisation.

 In the second part of the paper, we  define two families of Hopf algebras, denoted $\HCMlambda$, $\UCMlambda$,  parameterised by  $\l \in \field$. For $\l \neq 0$, the corresponding element of each family is isomorphic to the bicrossproduct $\HCM$ respectively $\UCM$.
  For $\l =0$ (the so-called classical limit) the Hopf algebra $\HCMlambda$ is commutative, and can be realised as functions  on the  semidirect product $\bR^2\lcross \Dplus$. 
 We construct a natural quotient Hopf algebra $\clh$ of $\HCMlambda$, which for $\l =0$ similarly corresponds  to  the coordinate algebra of the Heisenberg group. 
 For $\l \neq 0$ $\clh$ pairs with a  Hopf subalgebra $\ulh$ of $\UCMlambda$.
 By passing to an extended bicrossproduct $\ahcmdual \rlbicross  \FBplus_\l$ we give the correct classical limits of $\UCMlambda$ and $\ulh$.
   Finally we show that  $\clh$ and $\ulh$ are linked to  a local factorisation of the group $\SLtwoR$, in the same way $\HCM$ and $\UCM$ are linked to the factorisation of $\Diff^+ (\bR)$. 
    We remark that locally compact quantum groups (in the von Neumann algebra setting) similar  to $\clh$ and $\ulh$ were previously constructed by Vaes \cite{vaes}, linked to a factorisation of the continuous Heisenberg group rather than $\SLtwoR$.

Finally, a bicrossproduct coacts canonically on one of its factors (the Schr\"odinger coaction) 
hence a corollary of our results is that $\HCMlambda$ and $\clh$ coact canonically on $\Ulbplus$. 
The latter is  $\hhcm$  viewed as a noncommutative space, i.e. with scaling parameter $\l$ introduced in such a way as to be commutative when $\l=0$. 
This puts $\HCMlambda$ and $\clh$ in the same family as the coordinate algebras of the  Euclidean quantum  group of \cite{SM:thesis} and the $\kappa$-Poincar\'e quantum group \cite{lr,mr} coacting on algebras $U_\l ({{\bf b}}_{+}^{n})$ of various dimensions. 
In  such models  one is also interested in the covariant noncommutative differential geometry of the coordinate algebras of both the noncommutative space  and the coacting quantum group.
Thus in the final part of the paper we study low-dimensional covariant first order differential calculi over $U_\l (\bplus)$ and $\clh$.

\section{Preliminaries}

In this section we recall from \cite{SM_book} the construction of the (left-right) bicrossproduct Hopf algebra  $\A \lrbicross \H$ from 
    Hopf algebras $\A$ and $\H$, with $\A$ a left $\H$-module algebra and $\H$ a right $\A$-comodule coalgebra.
 For completeness we also give the definition of a  factorisation of a group $X$ into subgroups $G$ and $M$, and  the construction of a dual pair of finite-dimensional bicrossproduct Hopf algebras associated to a factorisation of a finite group 
 (this is not used directly in our constructions of infinite-dimensional bicrossproducts, but is an important part of the motivation).
We then  define the Hopf algebras $\ahcm$,  $\ahcmdual$, $\hhcm$, $\hhcmdual$  which we use to construct bicrossproducts.
 As shown by Figueroa and Gracia-Bondia  \cite{fgb}, $\ahcm$ is isomorphic to two other well-known Hopf algebras, the comeasuring Hopf algebra of the real line $\C$ \cite{SM_qbdg} and the Fa\`a di Bruno Hopf algebra $\F$.
  We use this to give a more convenient alternative  presentation (\ref{new_presentation_HCM}) of $\HCM$  using the generators $t_n$ of $\C$ instead of the $\d_n$.

\subsection{Bicrossproduct Hopf algebras}\label{section_bicrossproduct} 

Throughout this paper we work over a field $\field$ assumed to be of characteristic zero. 
For a Hopf algebra $\H$, we use the Sweedler notation $\D( x) = \sum x_\one \otimes x_\two$ for the coproduct.  We denote a right coaction $\D_R : \M \rto \M \otimes \H$ of $\H$ on a $\field$-vector space $\M$ by $\D_R ( m ) = \sum  m^\comid \otimes m^\coright$. Now let $\A$ be an algebra and $\C$ a coalgebra (over $\field$).

\begin{defn} $\A$ is a left $\H$-module algebra if there exists  a $k$-linear map 
$\tr : \H \otimes \A \rto \A$ such that $h \tr (ab) = \sum (h_\one \tr a)(h_\two \tr b)$ and  $h \tr 1 = \e(h) 1$, for all $h \in \H$ and $a$, $b \in \A$.
\end{defn}

\begin{defn} $\H$ is   a right $\C$-comodule coalgebra if there exists a right coaction $\D_R : \H \rto \H \otimes \C$ 
 such that for all $h \in \H$, 
 $$\sum \e( h^\comid ) h^\coright = \e(h) 1, \quad \sum {h^\comid}_{(1)} \otimes {h^\comid}_{(2)} \otimes h^\coright = \sum {h_{(1)}}^\comid \otimes {h_{(2)}}^\comid \otimes {h_{(1)}}^\coright \, {h_{(2)}}^\coright
$$
\end{defn}

\begin{defn}\label{dualpairing} We say that Hopf algebras $\H$, $\K$ are dually paired (in duality) if there exists a bilinear form $<.,.> : \H \times \K \rto \field$ such that
$$< a ,xy > = \sum < a_\one ,x > < a_\two , y>, \quad < ab ,x > = \sum < a , x_\one > < b , x_\two >$$
$$< S(a) , x> = < a , S(x) > , \quad < a,1> = \e(a) \quad < 1, x> = \e(x)$$
for all $a, b \in \H$, $x, y \in \K$. We say that the pairing is nondegenerate if for every nonzero $a \in \H$, there exists some $x \in \K$ so that $< a, x> \neq 0$, and for every nonzero $y \in \K$, there exists some $b \in \H$ so that $< b, y> \neq 0$. 
\end{defn}

If $\A$ and $\H$ are bialgebras, with $\H$ acting on $\A$, and $\A$ coacting on $\H$ (with action and coaction compatible in a suitable sense) then the bicrossproduct construction \cite{SM_hvn,SM_book} manufactures a larger bialgebra, the bicrossproduct of $\A$ and $\H$, containing both $\A$ and $\H$ as sub-bialgebras. If $\A$ and $\H$ are Hopf algebras, then so is the bicrossproduct.
Explicitly:

\begin{thm} \cite[Theorem 6.2.2]{SM_book}
\label{bicrossproduct_thm_1}
Let $\A$ and $\H$ be Hopf algebras, with $\A$ a left $\H$-module algebra, and $\H$ a right $\A$-comodule coalgebra, such that:
\begin{enumerate}
\item{
$\e( h \triangleright a) = \e(h) \e(a)$, $\D( h \triangleright a) = \sum {h_{(1)}}^{\comid} \triangleright a_{(1)} \otimes {h_{(1)}}^{\coright} ( h_{(2)} \triangleright a_{(2)} )$,
}
\item{\label{D_R_of_product}
$\D_R(1) = 1 \otimes 1$,
$\D_R (gh) = \sum {g_{(1)}}^{\comid}  h^{\comid} \otimes {g_{(1)}}^{\coright} ( g_{(2)} \triangleright h^{\coright} )$,
 }
 \item{
$\sum {h_{(2)}}^{\comid} \otimes ( h_{(1)} \triangleright a) {h_{(2)}}^{\coright} = \sum {h_{(1)}}^{\comid} \otimes {h_{(1)}}^{\coright} ( h_{(2)} \triangleright a )$.
}
\end{enumerate}
for all $a$, $b \in \A$, $g$, $h \in \H$. 
Then the vector space $\A \otimes \H$ can be given the structure of a Hopf algebra, the left-right bicrossproduct denoted  $\A \lrbicross \, \H$, via:
$$(a \otimes h)( b\otimes g) = \sum a (h_{(1)} \tr b) \otimes h_{(2)} g, \quad S(a \otimes h) = \sum ( 1 \otimes S h^{\comid} ) ( S ( a h^{\coright} ) \otimes 1)$$
\begin{equation}\label{bicross_antipode}\label{bicross_product_coproduct}
 S(a \otimes h) = \sum ( 1 \otimes S h^{\comid} ) ( S ( a h^{\coright} ) \otimes 1)
\end{equation}
\end{thm}

The left-right reversed result, constructing a Hopf algebra $\H \rlbicross \A$ from a right $\H$-module algebra $\A$ and  left $\A$-comodule coalgebra $\H$  is \cite[Theorem 6.2.3]{SM_book}.

\subsection{Group factorisations and finite-dimensional bicrossproducts}\label{group_factorisations}

A group $X$ is said to factorise into subgroups $G$ and $M$  if group multiplication gives a set bijection $G \times M \rightarrow X$. That is, given $x \in X$, there are unique $g \in G$, $m \in M$ such that $gm =x$. We write $X = G \bowtie M$.
Hence for any $m \in M$, $g \in G$ there exist unique  $g' \in G$, $m' \in M$ such that $mg = g' m'$. 
Writing $g' = m \triangleright g$, 
$m' = m \triangleleft g$, it is straightforward to check that this  defines a natural left action $\triangleright$ of $M$ on $G$, and a natural right action $\triangleleft$ of $G$ on $M$. 
For any group $\Gamma$ denote by $\field \Gamma$ the group algebra (over $\field$) of $\Gamma$, with coproduct $\D(g) = g \otimes g$, and (if $\Gamma$ is finite) by $\field [ \Gamma]$ the commutative algebra of $\field$-valued functions on $\Gamma$, with basis the delta-functions $\{ \d_g \}_{g \in \Gamma}$ and coproduct $\D ( \d_g ) = \sum_{xy = g} \d_x \otimes \d_y$.
 Then if a finite group $X$ factorises as  $X = G \bowtie M$, we can construct two dually-paired bicrossproducts:
  
 \begin{enumerate}
 \item $k[M] \lrbicross kG$.  We have a compatible left action and right coaction
  \begin{equation}\label{finite_group_one}
 g \tr \d_m := \d_{m \tl g^{-1}}, \quad g \mapsto \sum_{m \in M} (m \tr g) \otimes \d_m
 \end{equation}
   so we can form the bicrossproduct $k[M] \lrbicross kG$, with relations 
  $$g \d_m = \d_{(m \tl g^{-1} )} g, \quad \D( \d_m ) = \sum_{xy = m } \d_x \otimes \d_y, \quad \D(g) = \sum_{m \in M} (m \tr g) \otimes \d_m g$$
 
 \item $kM \rlbicross k[G]$. We have a compatible right action and left coaction
 $$\d_g \tl m := \d_{m^{-1} \tr g}, \quad m \mapsto \sum_{g \in G} \d_g \otimes (m \tl g)$$
  The bicrossproduct $kM \rlbicross k[G]$ has relations
 $$\d_g m = m \d_{m^{-1} \tr g}, \quad \D( \d_g) = \sum_{xy=g} \d_x \otimes \d_y, \quad \D(m) = \sum_{g \in G} m \d_x \otimes (m \tl g)$$
 \end{enumerate}
 
 When $G$ and $M$ are infinite, both these  constructions fail. The challenge for factorisations of infinite groups is to construct bicrossproduct Hopf algebras  analogous to those above.
  In Section \ref{section_HCM_is_bicrossproduct} we solve this problem for the group  factorisation $\Diff^+ (\bR) = \Bplus \bowtie \Dplus$.

\subsection{The commutative Hopf subalgebra $\ahcm$}
\label{comm_Hopf_subalg}

We define $\ahcm$  to be  the unital commutative  subalgebra of $\HCM$ generated by $\{ \d_n \; : \; n =1 ,2 , \ldots \}$, which, as shown in \cite{cm}, is a Hopf subalgebra of $\HCM$ with 
$\D( \d_n) = \sum_{k=1}^n D_{n,k} \otimes \d_k$
for some $D_{n,k} \in \ahcm$. It was shown in \cite{fgb} that $\ahcm$ is isomorphic 
 to both the comeasuring Hopf algebra $\C$ of the real line  and the Fa\`a di Bruno Hopf algebra $\F$, whose definitions we now recall.

\begin{defn}\label{comeasuring} \cite{SM_qbdg} The comeasuring Hopf algebra $\C$ of the real line is the commutative Hopf algebra over $\field$ generated by indeterminates $\{ \; t_n \; : \; n=1,2, \ldots \; \}$ with $t_1 =1$,  counit $\e ( t_n ) = \d_{n,1}$, and coproduct  
\begin{equation}
\label{Delta_t_n_first}
\D( t_n ) = 
 \sum_{k=1}^n \; ( \sum_{i_1 + \ldots + i_k = n} t_{i_1} \ldots t_{i_k}) \otimes t_k
\end{equation}
Adapting results of  \cite{fgb}, the  antipode on $\C$ is given by
\begin{equation}
\label{antipode_t_n}
S( t_{n+1}) = \sum_{{\bf c} \in \S} (-1)^{n - c_1} {\frac{(2n- c_1)! c_1 !}{(n+1)!}} \; {\frac{ t_1^{c_1} t_2^{c_2} \ldots t_{n+1}^{c_{n+1}}}{ c_1 ! c_2 ! \ldots c_{n+1} !}}
\end{equation}
where $\S = \{ (c_1, \ldots , c_{n+1} ) \; : \; \sum_{j=1}^{n+1} c_j = n, \;  \sum_{j=1}^{n+1} j c_j = 2n \; \}$.
\end{defn}

  If we rewrite Definition \ref{comeasuring} in terms of  generators $a_n =  n! \; t_n$, this gives  the usual presentation of the Fa\`a di Bruno Hopf algebra $\F$.

\begin{prop}\label{isomorphism_ahcm} \cite{fgb}
$\ahcm$, the Fa\`a di Bruno Hopf algebra $\F$ and the comeasuring Hopf algebra of the real line $\C$ are all  isomorphic, via
$$\d_n \mapsto n! \, \sum_{{\bf c} \in \S} (-1)^{n-c_1 } {\frac{ (n - c_1 )! }{ c_2 ! \ldots c_{n+1} !}} (t_1)^{c_1} (2 t_2 )^{c_2} \ldots ( (n+1) t_{n+1} )^{c_{n+1}}$$
where $\S = \{ (c_1, \ldots , c_{n+1} ) \;:\; \sum_{j=1}^{n+1} c_j =  n+1, \; \sum_{j=1}^{n+1} j c_j  = 2n+1 \}$, and 
$$(n+1) t_{n+1} \mapsto  \sum_{c_1 + 2 c_2 + \ldots + n c_n =n} {\frac{\d_1^{c_1} \ldots \d_n^{c_n}}{ c_1 ! \ldots c_n ! (1!)^{c_1} \ldots (n!)^{c_n}}}$$
\end{prop}

In the sequel, we will use the presentation of $\ahcm$ as the commutative Hopf algebra with generators $\{ t_n \}_{n \geq 1}$, and coproduct and antipode given by (\ref{Delta_t_n_first},\ref{antipode_t_n}).

\begin{lemma}\label{ahcm_grading}
 $\ahcm$ is an $\bN$-graded Hopf algebra, via the grading defined on monomials by
$| t_{n_1} \ldots t_{n_A} | = n_1 + \ldots + n_A - A$.
\end{lemma}
\begin{prf} Let $\ahcm_N$ be the linear span of monomials ${\bf t}= t_{n_1} \ldots t_{n_A}$ with $|{\bf t}| = N$.
 Then  $\ahcm_m \ahcm_n \subseteq \ahcm_{m+n}$ for all $m$, $n$. 
From (\ref{antipode_t_n})  $S( \ahcm_N ) \subseteq \ahcm_N$.
It also follows from 
 (\ref{Delta_t_n_first}) that 
$\D( \ahcm_N ) \subseteq \oplus_{n=0}^N \; \ahcm_n \oplus \ahcm_{N-n}$. 
\end{prf}

As a corollary of Proposition \ref{isomorphism_ahcm}, we have:

\begin{cor} In terms of the $t_n$, the presentation (\ref{HCM_coproduct}) of $\HCM$ becomes:
\begin{eqnarray}
&&[Y, X] = X, \quad [ X, t_n ] = (n+1) t_{n+1} - 2 t_2 t_n, \quad [Y, t_n ] = (n-1) t_n\nonumber\\
\label{HCM_X_Y_t_n}
&&\D( X) = X \otimes 1 + 1 \otimes X +  Y \otimes 2 t_2, \quad \D( Y) = Y \otimes 1 + 1 \otimes Y, \quad \e( t_n ) = \d_{n,1}\nonumber\\
\label{new_presentation_HCM}
&&S(X) = - X + 2 Y t_2, \quad S(Y) = -Y,\quad \e(X) = 0 = \e(Y)
\end{eqnarray}
with $\D( t_n)$, $S(t_n)$ given by (\ref{Delta_t_n_first}, \ref{antipode_t_n}).
\end{cor} 

 Taking $\field = \bR$ or $\bC$, the generators $\d_n$, $t_n$ can be realised as  functions on $\Dplus$: 
  \begin{equation}
\label{d_n_t_n}
\d_n (f ) = [ \log f' ]^{(n)} (0), \quad t_n (f) = {\frac{1}{n!}}f^{(n)}(0), \quad f \in \Dplus
\end{equation}

\subsection{The Hopf algebra $\ahcmdual$}

\begin{defn} We define ${\bf d}_+$ to be  the Lie algebra (over $\field$) with countably many generators $\{ z_n \}_{n \in \bN}$   and relations $[ z_m , z_n ] = (n-m) z_{m+n-1}$ for all $m$, $n$. Define 
 $\dplus$ to be the Lie subalgebra  generated by  $\{ z_n \}_{n \geq 2}$.
 Then $U(\dplus)$ is the  universal enveloping algebra of $\dplus$ with canonical Hopf structure.
\end{defn}

Note that ${\bf d}_+$  has a natural representation as differential operators $z_n = x^n {\frac{d}{dx}}$ acting on the unital algebra $k[x]$ of polynomials in a single indeterminate $x$.

\begin{lemma}\label{Udplus_grading} $U(\dplus)$ is an $\bN$-graded Hopf algebra, via the grading defined on monomials by $|1|=0$, $| z_{m_1} \ldots z_{m_p} | = m_1 + \ldots + m_p - p$.
\end{lemma}
\begin{prf} Denote by ${U(\dplus)}_N$ the linear span of monomials of degree $N$. 
Since $| z_m z_n | = m+n -2 = |z_n z_m | =  | z_{m+n-1} |$, then ${U(\dplus)}_m {U(\dplus)}_n \subseteq {U(\dplus)}_{m+n}$ for all $m$, $n \in \bN$.
 Further, as $\D( z_n ) = z_n \otimes 1 + 1 \otimes z_n$, it follows that 
 $\D( {U(\dplus)}_N ) \subseteq \oplus_{n=0}^N \; {U(\dplus)}_n \otimes {U(\dplus)}_{N-n}$.
  Finally $S( {U(\dplus)}_N ) \subseteq {U(\dplus)}_N$.
\end{prf}

\begin{prop}\label{z_m_pairs_with_t_n} \label{z_m_pairs_with_d_n}
 There is a nondegenerate dual pairing (in the sense of Definition \ref{dualpairing}) of the Hopf algebras $\ahcmdual$ and $\ahcm$, defined on generators by 
 \begin{equation}
 \label{pairing_A_Adual}
< z_m , t_n > =  \d_{m, n} \quad \quad \forall \; m \geq 2, \; n \geq 1
\end{equation}
equivalently by $< z_m , \d_n > = m! \; \d_{m,n+1}$. The pairing satisfies
\begin{eqnarray}
&< z_{m_1}  \ldots z_{m_p} , t_n > &=
\left\{
	\begin{array}{ll}
	\Pi_{j=1}^{p-1} ( n+j-1 - \sum_{l=1}^j m_l ) 
	 &: \sum_{j=1}^{p} m_j = n+p-1\\
	0 &: \otherwise
	\end{array}\right.\nonumber\\
\label{pairing_satisfies}
&< z_m , t_{n_1} \ldots t_{n_A} > &=
  \left\{
	\begin{array}{ll}
	1 &: \{ n_1 , \ldots , n_A \} = \{ m,1, \ldots , 1\} \; \assets\\
	0 &: \otherwise
	\end{array}\right.
\end{eqnarray}
\end{prop}
\begin{prf} Assuming the pairing is well-defined, the identities follow by a straightforward induction. 
For example, 
$$< z_m , t_{n_1} \ldots t_{n_A} > = < \D^{A-1}( z_m ) , t_{n_1} \otimes \ldots \otimes t_{n_A} > = \sum_{l=1}^A \e( t_{n_1} ) \ldots < z_m , t_{n_l}> \ldots \e( t_{n_A} )$$
To check well-defined, as
$< z_m z_n , t_{n_1} \ldots t_{n_A} > =$
$< \D^{A-1} ( z_m z_n ) , t_{n_1} \otimes \ldots \otimes t_{n_A}>$ 
 it suffices to check $< z_m z_n , t_p >$. 
By the above, $< z_m z_n , t_p > = n \d_{m+n ,p+1}$. 
So $< z_m z_n - z_n z_m , t_p > = (n-m) \d_{m+n-1, p} = (n-m) < z_{m+n-1}, t_p >$.

 It follows from (\ref{antipode_t_n},\ref{pairing_satisfies}) that $<S(z_m), t_{n}> =  < z_m , S(t_{n}) >$ for all $m$, $n$, so 
$<S(z_m), {\bf t}> =$ $< z_m , S({\bf t}) >$ for ${\bf t} = t_{n_1} \ldots t_{n_A}$ hence for all ${\bf t} \in \hhcmdual$. Then
 $$< S( z_{m_1} \ldots z_{m_p} ) , {\bf t} > = < S( z_{m_p}) \otimes \ldots \otimes S( z_{m_1}), \D^{p-1} ( {\bf t} )>$$
  $$= < z_{m_1} \otimes \ldots \otimes z_{m_p} , \D^{p-1} ( S({\bf t}))> = < z_{m_1} \ldots z_{m_p} , S( {\bf t})>$$

To show nondegeneracy, we need the $\bN$-gradings of  Lemmas \ref{ahcm_grading} and  \ref{Udplus_grading}.

\begin{lemma}\label{nondeg_lemma} For ${\bf z} = z_{m_1}^{a_1} \ldots z_{m_p}^{a_p}$ and ${\bf t} = t_{n_1} \ldots t_{n_A}$, with $2 \leq m_1 < m_2 < \ldots < m_p$, $2 \leq n_1 \leq n_2 \leq \ldots \leq n_A$, and $a_1, \ldots a_p , A \geq 1$, then:
\begin{enumerate}
\item\label{ztone} $< {\bf z}, {\bf t} > = 0$ if $A > a_1 + \ldots + a_p$.
\item\label{zttwo} $< {\bf z}, {\bf t} > = 0$ unless $| {\bf z}| = |{\bf t}|$, i.e. unless $a_1 m_1 + \ldots + a_p m_p = n_1 + \ldots + n_A$.
\item\label{ztthree} If $A = a_1 + \ldots + a_p$,  $< {\bf z}, {\bf t} > = a_1 ! \ldots a_p ! \d_{m_1 , n_1} \ldots \d_{m_1 , n_{a_1}} \d_{m_2 , n_{a_1 + 1}} \ldots \d_{m_p , n_A}$.
\end{enumerate}
\end{lemma}
\begin{prf} For part \ref{ztone}, using $< {\bf z} , {\bf t} > = < \D^{A-1} ( {\bf z} ), t_{n_1} \otimes \ldots \otimes t_{n_A} >$, if $A > a_1 + \ldots + a_p$ then every term in $\D^{A-1} ( {\bf z} )$ will contain at least one component $ - \otimes 1 \otimes -$, which pairs to zero with the corresponding $t_{n_k}$. Hence $< {\bf z} , {\bf t} > =0$.
As $\ahcm$ and  $U(\dplus)$ are both $\bN$-graded, part \ref{zttwo} follows using  (\ref{pairing_satisfies}).
We prove part \ref{ztthree} by induction. It holds for $p=A=1$. Suppose it holds for ${\bf z}$, ${\bf t}$. Then for $m_p \leq m_{p+1}$, $n_A \leq n_{A+1}$, we have
$$< {\bf z} \, z_{m_{p+1}} , {\bf t} \, t_{n_{A+1}}> = < {\bf z} \otimes z_{m_{p+1}} , \D( {\bf t} \, t_{n_{A+1}}) >
 = < {\bf z} \otimes z_{m_{p+1}} , \sum \ldots \otimes t_{i_1} \ldots t_{i_{A+1}} >$$
 \begin{equation}\label{usethisagain}
 = \sum_{k=1}^A < {\bf z} , t_{n_1} \ldots \hat{t_{n_k}} \ldots t_{n_{A+1}} > \d_{m_{p+1}, n_k } + < {\bf z} , {\bf t} >  \d_{m_{p+1}, n_{A+1}}
 \end{equation}
 using (\ref{pairing_satisfies}).
  The first $A$ terms contain $\d_{m_p, n_{A+1}} \d_{m_{p+1}, n_k}$.   If $m_p < m_{p+1}$ this is zero, so 
   $< {\bf z} \, z_{m_{p+1}} , {\bf t} \, t_{n_{A+1}}> = < {\bf z}, {\bf t}> \d_{m_{p+1}, n_{A+1}}$.
   If $m_p = m_{p+1}$, then (\ref{usethisagain}) becomes
   $$\sum_{k=A - a_p +1}^A < {\bf z}, t_{n_1} \ldots \hat{t_{n_k}} \ldots t_{n_{A+1}} > \d_{m_{p}, n_k } + < {\bf z} , {\bf t} >  \d_{m_{p}, n_{A+1}} = ( a_p +1 ) <{\bf z} , {\bf t} > \d_{m_{p}, n_{A+1}}$$
    which completes the inductive step.
   \end{prf}

We can now prove nondegeneracy. Given ${\bf z} = \sum \l( {\bf m}, {\bf a}) z_{m_1}^{a_1} \ldots z_{m_p}^{a_p}$ such that  $< {\bf z}, {\bf t}> =0$ for all ${\bf t} \in \ahcm$, 
by Lemma \ref{nondeg_lemma}, part \ref{zttwo}  we can restrict to ${\bf z} \in {U(\dplus)}_N$ for some $N$. 
 Define $A = \max \{ a_1 + \ldots + a_p \}$ taken over monomials occurring in ${\bf z}$. 
  By Lemma \ref{nondeg_lemma}, part \ref{ztone},
 $$< {\bf z}, {\bf t} > = < \sum_{ {\bf m}, {\bf a} \; : \; a_1 + \ldots + a_p = A} \l ({\bf m} , {\bf a} ) z_{m_1}^{a_1} \ldots z_{m_p}^{a_p} , t_{n_1} \ldots t_{n_A} >$$
  So  by  Lemma \ref{nondeg_lemma}, part \ref{ztthree},  choosing ${\bf t}$ appropriately gives $\l({\bf m},{\bf a}) = 0$ for all such ${\bf m}$, ${\bf a}$. 
 Nondegeneracy follows,  completing the proof of Proposition \ref{z_m_pairs_with_t_n}.
\end{prf}

\subsection{The  Hopf algebras $\hhcm$ and  $\hhcmdual$}
\label{hhcm_dual}

The  group $\Bplus$  was defined in    (\ref{defn_Bplus}).
Let $\bplus$ be the Lie algebra (over $k$) generated by $X$, $Y$ satisfying the  relation $[Y, X] = X$.
We now define a commutative  Hopf algebra  $\hhcmdual$ and a  nondegenerate dual pairing  of $\hhcmdual$ and $\hhcm$.

\begin{defn}  $\hhcmdual$ is the commutative Hopf algebra (over $\field$)   generated by 
 elements $\a^{\pm 1}$, $\b$ satisfying 
 \begin{eqnarray}
&&\D(\a) = \a \otimes \a, \quad
\D(\b) = \a \otimes \b + \b \otimes 1\nonumber\\
&&\e( \a) =1, \quad \e(\b)=0,\quad S(\a) = \a^{-1}, \quad S( \b ) = - \a^{-1} \b 
\nonumber
\end{eqnarray}
\end{defn}

\begin{lemma} There is a unique nondegenerate dual pairing of the 
 Hopf algebras $\hhcm$  and $\hhcmdual$, defined on generators by 
  \begin{eqnarray}
  &&< X, \a > = 0, \quad < X, \b> = 1= <Y, \a >, \quad < Y, \b > =0\nonumber
  \end{eqnarray}
 This satisfies 
\begin{equation}\label{pairing_hhcm_FBplus}
\label{pairing_X^j_Y^k_with_a^p_b^q}
< X^j Y^k , \a^t  \b^r > =  j! \d_{j,r}  t^k  \quad \forall \; j,k,r \in \bN, \; t \in \bZ
\end{equation}
where we use the convention  $X^0 =1=Y^0 = \a^0 = \b^0$, $0! =1$, and $0^0 =1$.
\end{lemma}
\begin{prf} This is more  straightforward than the proof of Proposition \ref{z_m_pairs_with_t_n}, and also well known. We give the details for completeness.
For $t \geq 1$, $< Y, \a^t > = \sum_{i=1}^t < 1 \otimes \ldots Y \ldots \otimes 1 , \a^{\otimes t} > = t$. 
We also have $< Y, \a^{-1}> = < Y, S(\a) > = < S(Y) , \a > = - < Y, \a > = -1$, and in fact $< Y, \a^t > = t$ for all $t \in \bZ$. 
Then  $< Y^k , \a^t > = < Y^{\otimes k} , \D^{k-1} ( \a^t)>$ $= < Y, \a^t >^k = t^k$. 
Suitably interpreted, this holds also for $k=0$.
In the same way, $< Y^k , \b^r > = \d_{k,0} \d_{r,0}$, for all $k$, $r \geq 0$. So
$< Y^k , \a^t \b^r > = \d_{r,0} t^k$. 
Further, $< X^j , \a^t > = \d_{j,0}$  and $< X^j \b^r > = j! \d_{j,r}$
for all $j$, $r, t \geq 0$, hence
$< X^j , \a^t \b^r > = j! \,  \d_{j,r}$.
Using  
$< X^j Y^k , \a^t \b^r > =$ $ < X^j \otimes Y^k , \D(  \a^t  \b^r )>$
the result follows. 
  It is also straightforward to check that $< S( X^j Y^k ), \a^t  \b^r > = < X^j Y^k , S( \a^t  \b^r )>$.

To check well-defined, we have $<XY, \a^t \b^r > = \d_{1,r} t$, and
\begin{eqnarray}
&< YX, \a^t \b^r > &= < Y \otimes X, ( \a^t \otimes \a^t ) ( \sum_{s=0}^r {\tiny{\rchooses}} \a^s \b^{r-s} \otimes \b^s ) >\nonumber\\
&&= \sum_{s=0}^r {\tiny{\rchooses}} \d_{0,r-s} (s+t) \d_{1,s} = \d_{1,r} (t+1)\nonumber
\end{eqnarray}
Hence $< YX - XY, \a^t \b^r > = \d_{1,r} = < X, \a^t \b^r >$.
 Finally, given (\ref{pairing_X^j_Y^k_with_a^p_b^q}), nondegeneracy of the pairing is immediate.
\end{prf}

\section{Bicrossproduct structure of the Connes-Moscovici Hopf algebra}
\label{section_HCM_is_bicrossproduct}

We prove that $\ahcm$ can be given the structure of a left $\hhcm$-module algebra, and $\hhcm$ the structure of a right $\ahcm$-comodule coalgebra, with action and coaction compatible in the sense of Theorem \ref{bicrossproduct_thm_1}. This enables the construction of a bicrossproduct Hopf algebra $\hcm$, which we prove is isomorphic to the Connes-Moscovici Hopf algebra $\HCM$.
  We also construct a second bicrossproduct $\UCM := \hcmdual$, equipped with a nondegenerate dual pairing with $\HCM$.
 We explain how these bicrossproducts are linked to the  factorisation of the 
 group $\Diff^+ (\bR)$ 
   into the subgroups $\Bplus$ and $\Dplus$. 
   The factorisation argument is not part of our proof, but rather serves as motivation.

 \subsection{The bicrossproduct $\hcm$}
\label{section_hcm}

 \begin{lemma}
  \label{lemma_left_action_H_on_A}
  $\ahcm$ is a left $\hhcm$-module algebra via the action defined by
\begin{equation}\label{left_action_H_on_A}
X \tr t_n = (n+1) t_{n+1} - 2 t_2 t_n, \quad Y \tr t_n = (n-1) t_n
\end{equation} 
equivalently defined by $X \tr \d_n  =  \d_{n+1}$, $Y \tr \d_n  =  n \d_n$.
  \end{lemma}
 \begin{prf} As $\ahcm$ is commutative and $\hhcm$ cocommutative, it  is easy to check that  (\ref{left_action_H_on_A})
  extends to a well-defined left action of $\hhcm$ on $\ahcm$. For example, 
 \begin{eqnarray}
 &Y \tr (X \tr t_n ) &= Y \tr [ (n+1) t_{n+1}  -2 t_2 t_n ] =   (n^2 + n) t_{n+1} - 2n t_2 t_n,\nonumber\\
 &X \tr (Y \tr t_n) &= (n-1) X \tr t_n = (n^2 -n) t_{n+1} - 2(n-1) t_2 t_n\nonumber
 \end{eqnarray}
Hence $(YX - XY) \tr t_n = (n+1) t_{n+1} - 2 t_2 t_n = X \tr t_n$. 
The action on the $\d_n$ follows from Proposition \ref{isomorphism_ahcm}.
\end{prf}

\begin{lemma}
\label{lemma_right_coaction_A_on_H}
$\hhcm$ is a right $\ahcm$-comodule coalgebra, via the coaction defined on generators by 
\begin{equation}
\label{right_coaction_A_on_H}
\D_R (X) = X \otimes 1 + Y \otimes 2 t_2, \quad \D_R (Y) = Y \otimes 1
\end{equation}
and extended by 
$\D_R (gh) = \sum {g_{(1)}}^{\comid}  h^{\comid} \otimes {g_{(1)}}^{\coright} ( g_{(2)} \triangleright h^{\coright} )$.
\end{lemma}
\begin{prf} We check that these formulae define a coaction. 
Suppose $g \in \hhcm$ satisfies $(\id \otimes \D) \D_R (g) = (\D_R \otimes \id) \D_R (g)$ (this holds for $g = X, Y$). Then 
\begin{eqnarray}
&&\D_R (gY) = \sum {g_{(1)}}^{\comid}  Y^{\comid} \otimes {g_{(1)}}^{\coright} ( g_{(2)} \triangleright Y^{\coright} )
 = \sum {g}^{\comid}  Y \otimes {g}^{\coright}\nonumber\\
&&\Rightarrow (\D_R \otimes \id) \D_R (gY) = \sum {g}^{\comid \comid}  Y \otimes {g}^{(1) \coright} \otimes g^{(2)}\nonumber\\
&&= \sum g^\comid Y \otimes {g^\coright}_\one \otimes {g^\coright}_\two =
(\id \otimes \D) \D_R (gY)\nonumber
\end{eqnarray}
and in the same way $(\D_R \otimes \id) \D_R (Xg) = (\id \otimes \D) \D_R (Xg)$.
We check that $\D_R$ is well-defined. We have
\begin{eqnarray}
&&\D_R (YX) = \sum {Y_{(1)}}^{\comid}  X^{\comid} \otimes {Y_{(1)}}^{\coright} ( Y_{(2)} \triangleright X^{\coright} )= YX \otimes 1 + (Y^2 + Y)\otimes 2 t_2\nonumber\\
&&\D_R (XY) = XY \otimes 1 +  Y^2 \otimes 2 t_2\nonumber
\end{eqnarray}
Hence $\D_R ( YX - XY ) = (YX - XY ) \otimes 1 + Y \otimes 2 t_2 = \D_R (X)$.
\end{prf}

\begin{thm}\label{compatible}
 The left action (\ref{left_action_H_on_A}) and right coaction (\ref{right_coaction_A_on_H}) are compatible  in the sense of Theorem \ref{bicrossproduct_thm_1}.
\end{thm}
\begin{prf} We check conditions 1-3 of Theorem \ref{bicrossproduct_thm_1}.
 For 1, to show $\e(h \tr a) = \e(h) \e(a)$, it is enough to show $\e( X \tr a) = 0 = \e( Y \tr a)$ for all $a$.
  Using the $\bN$-grading of $\ahcm$ (Lemma \ref{ahcm_grading}) we see that 
  $X \tr \ahcm_N \subseteq \ahcm_{N+1}$, $Y \tr \ahcm_N \subseteq \ahcm_N$, and $\e( \ahcm_N ) =0$ for $N \geq 1$.
 We also need to check that 
\begin{equation}
\label{Delta_h_tr_a}
\D( h \triangleright a) = \sum {h_{(1)}}^{\comid} \triangleright a_{(1)} \otimes {h_{(1)}}^{\coright} ( h_{(2)} \triangleright a_{(2)} )
\end{equation}
for all $h \in \hhcm$, $a \in \ahcm$. 
Suppose for fixed $h$, $k$  (\ref{Delta_h_tr_a}) holds for all $a$. Then
$$\D( hk \tr a) = \D( h \tr (k \tr a)) =\sum {h_{(1)}}^\comid \tr ( k \tr a)_{(1)} \otimes {h_{(1)}}^\coright ( h_{(2)} \tr (k \tr a)_{(2)} )$$
$$= \sum ({h_{(1)}}^\comid {k_{(1)}}^\comid ) \tr a_{(1)} \otimes {h_{(1)}}^\coright ( h_{(2)(1)} \tr {k_{(1)}}^\coright ) ( h_{(2)(2)} k_{(2)} \tr a_{(2)} )$$
Now, $(\D_R \otimes \id) \D(hk) = \sum { (hk)_{(1)}}^\comid \otimes {(hk)_{(1)}}^\coright \otimes (hk)_{(2)} $
\begin{eqnarray}
&& = \sum {h_{(1)(1)}}^\comid {k_{(1)}}^\comid \otimes {h_{(1)(1)}}^\coright ( h_{(1)(2)} \tr {k_{(1)}}^\coright) \otimes h_{(2)} k_{(2)}
\nonumber\\
&& = \sum {h_{(1)}}^\comid {k_{(1)}}^\comid \otimes {h_{(1)}}^\coright ( h_{(2)(1)} \tr {k_{(1)}}^\coright ) \otimes  h_{(2)(2)} k_{(2)} 
\nonumber\\
&&\Rightarrow \quad \D( hk \tr a) = \sum {(hk)_{(1)}}^\comid \tr a_{(1)} \otimes  {(hk)_{(1)}}^\coright ( (hk)_{(2)} \tr a_{(2)} )\nonumber
\end{eqnarray}
So it will be enough to check $h = X$, $Y$ only. It is straightforward to check that if (\ref{Delta_h_tr_a}) holds for $\D( X \tr a)$, $\D(X \tr b)$, then it holds for $\D( X \tr ab)$, and similarly for $Y$. Hence we need only check $X \tr t_n$, $Y \tr t_n$. 
Now, $\D( Y \tr t_n ) = (n-1) \D(t_n)$, whereas the right-hand side of (\ref{Delta_h_tr_a}) is
\begin{eqnarray}
&&\sum {Y_{(1)}}^{\comid} \triangleright {t_n}_{(1)} \otimes {Y_{(1)}}^{\coright} ( Y_{(2)} \triangleright {t_n}_{(2)} ) = \sum [ Y \tr ({t_n})_{(1)} \otimes ({t_n})_{(2)} + ({t_n})_{(1)} \otimes Y \tr ({t_n})_{(2)}]\nonumber\\
&=& \sum_{k=1}^n \sum_{i_1 + \ldots + i_k =n} [ Y \tr ( t_{i_1} \ldots t_{i_k} ) \otimes t_k + t_{i_1} \ldots t_{i_k}  \otimes Y \tr t_k ]\nonumber\\
&=&  \sum_{k=1}^n \sum_{i_1 + \ldots + i_k =n} ( i_1 + \ldots + i_k -1) t_{i_1} \ldots t_{i_k}  \otimes t_k = (n-1) \D(t_n) = \D( Y \triangleright t_n)\nonumber
\end{eqnarray}
In the same way, we check that (\ref{Delta_h_tr_a}) holds for $\D(X \tr t_2)$, with the general case $\D(X \tr t_n )$ following by induction. Hence condition 1 of Theorem \ref{bicrossproduct_thm_1} holds.

Condition 2 is automatic from the definition of $\D_R$. 
Finally, since $\ahcm$ is commutative and $\hhcm$ cocommutative, condition 3 is immediate.
\end{prf}

All conditions of Theorem \ref{bicrossproduct_thm_1} hold, so we can construct the left-right bicrossproduct Hopf algebra $\hcm$, which we can think of as  an analogue of the finite-dimensional bicrossproduct $k[M] \lrbicross kG$ defined in Section \ref{group_factorisations}.

\begin{thm}\label{HCM_is_bicrossproduct}
 $\HCM$ and $\hcm$ are isomorphic Hopf algebras.
\end{thm}
\begin{prf} Using (\ref{bicross_product_coproduct}), it follows that $\hcm$ has generators $X$, $Y$, $t_n$  with relations  coinciding exactly with the presentation (\ref{new_presentation_HCM}) of $\HCM$.
\end{prf}

Finally, we remark that the original codimension one Connes-Moscovici Hopf algebra $\HCMoriginal$, which differs from  $\HCM$  only in that $Y \otimes \d_1$ is replaced by $\d_1 \otimes Y$ in (\ref{HCM_coproduct}), is isomorphic to  a right-left bicrossproduct $\hhcm \rlbicross \ahcm$ which is also linked to the factorisation of $\Diff^+ (\bR)$, but in a less natural way than $\HCM$. For  completeness, this is outlined (without proofs) in Section \ref{section_HCM_as_bicross}.

\subsection{Relation to group factorisation I}
\label{relation_to_gp_fact}

We motivate the above bicrossproduct constructions  using a factorisation of the 
 group $\Diff^+ (\bR)$ of orientation preserving diffeomorphisms of the real line (\ref{defn_Diffplus}).
  As shown in \cite{cm},  the factorisation $\Diff^{+} ( \bR )= \Bplus \bowtie \Dplus$ is as follows.
 Given $\varphi \in \Diff^+ (\bR)$, we have 
 $\varphi =  (a,b) \circ \phi$ for unique $(a,b) \in \Bplus$, $\phi \in \Dplus$, with 
 \begin{equation}
 \label{D=B+_bowtie_D0}
 (a,b) = ( \varphi'(0), \varphi(0)), \quad \phi(x) ={ \frac{ \varphi(x) - \varphi(0)}{\varphi' (0)}} \quad \forall \; x \in \bR
 \end{equation}
  Since $(\phi \circ (a,b))(x) = \phi(ax+b)$, the corresponding left action of  $\Dplus$ on $\Bplus$ and right action of $\Bplus$ on $\Dplus$ are given by 
 \begin{equation}
 \label{left_action_D0_on_B+}
 \label{right_action_B+_on_D0}
 \phi \tr (a,b) = ( a \phi'(b), \phi(b)), \quad (\phi \tl (a,b))(x) = {\frac{\phi (ax+b) - \phi(b)}{a \phi'(b)}}
 \end{equation}

We identify $X$, $Y \in \bplus$ with the matrices
${\tiny{\left(
\begin{array}{cc}
0 & 1\cr 0 & 0 \end{array}\right)}}$, 
${\tiny{\left(
\begin{array}{cc}
1 & 0\cr 0 & 0 \end{array}\right)}}$.
By slight abuse of notation, for any $s \in \bR$ denote by $e^{sX}$, $e^{sY}$ the elements $(1,s)$, $(e^s , 0)$ of $\Bplus$.

 To understand the origin of (\ref{left_action_H_on_A}), consider the factorisation (\ref{D=B+_bowtie_D0}) of $\Diff^{+} ( \bR )$. 
 For any function $\xi : \Dplus \rto \field$, using (\ref{left_action_D0_on_B+}) we define a  left action of $\Bplus$ via 
 $( (a,b) \tr \xi )(\phi) = \xi ( \phi \tl (a,b) )$, and (taking $k = \bR$ or $\bC$) by differentiation a left action of $\hhcm$ on $\ahcm$.
 So $(\phi \tl e^{sX}) (x) = {\frac{ \phi(x+s) - \phi(s)}{ \phi'(s)}}$, then (\ref{d_n_t_n}) gives
 \begin{eqnarray}
 &t_n ( \phi \tl e^{sX} ) &= {\frac{\phi^{(n)} (s)} { n! \phi' (s) }}= {\frac{[ \phi^{(n)}(0) + s  \phi^{(n+1)}(0) + O( s^2 ) ]}{ n! [ \phi' (0) + s \phi'' (0) + O (s^2 ) ] }}\nonumber\\
 &&= {\frac{1}{n!}} [ \phi^{(n)} (0) - s \phi'' (0) \phi^{(n)} (0) + s \phi^{(n+1)} (0)] + O (s^2 )\nonumber\\
 &\Rightarrow (X \tr t_n )(\phi)  &= {{\frac{d}{ds}} |}_{s=0} t_n ( \phi \tl e^{sX} ) =  {\frac{1}{n!}} [  -  \phi'' (0) \phi^{(n)} (0) +  \phi^{(n+1)} (0)]\nonumber\\
 &&= [ (n+1) t_{n+1}  - 2 t_2 t_n ] (\phi)\nonumber
 \end{eqnarray}
 giving $X \tr t_n = (n+1) t_{n+1}  - 2 t_2 t_n$. Similarly, 
 $\d_n ( \phi \tl e^{sX} ) =g^{(n)} (0)$, where 
 $$g(x) = \log (\phi \tl e^{sX} )' (x) = \log \phi' (x+s) - \log \phi' (s)$$
 Hence $g^{(n)} (x) = h^{(n)} (x+s)$, where $h(x) =  \log \phi'  (x)$, so $h^{(n)} (0) = \d_n (\phi)$ for all $n \geq 1$.
 So $g^{(n)} (0) = h^{(n)} (s) = h^{(n)} (0) + s h^{(n+1)} (0) + O( s^2 )$. Thus, 
 $$(X \tr \d_n )(\phi) =  {{\frac{d}{ds}}|}_{s=0} \d_n ( \phi \tl e^{sX} ) = h^{(n+1)} (0) = \d_{n+1} ( \phi)$$
So $X \tr \d_n = \d_{n+1}$. The formulae for $Y$ follow similarly. 
So using the group factorisation we  recover (\ref{left_action_H_on_A}), which we already showed to be a left action.

Next we explain how the coaction (\ref{right_coaction_A_on_H}) can be recovered from the  factorisation (\ref{D=B+_bowtie_D0}).
 For any group factorisation $X = G \bowtie M$ we define a $\field$-linear map 
 $$\tilde{\D}_R : \field G \rto  \field[M, \field G], \quad \tilde{\D}_R (g) (m) = m \tr g$$
 where $\field[M, \field G]$ is the $\field$-vector space of maps $M \rto \field G$. 
If $X$ is finite then $\field [ M , \field G ] \cong \field G \otimes \field [M]$ as vector spaces, and 
  $\tilde{\D}_R$ is the right coaction (\ref{finite_group_one}).
  In our situation, taking $k = \bR$ then for any $s \in \bR$, $\phi \in \Dplus$, 
  $$\tilde{\D}_R ( e^{sX} ) (\phi) = \phi \tr e^{sX} = \phi \tr (1,s) = ( \phi'(s), \phi(s))$$
   We induce a linear map  $\D_R : \hhcm \rto \bR [ \Dplus, \hhcm ] $ by differentiation:
  $$\D_R (X) (\phi) = {{\frac{d}{ds}}|}_{s=0} \tilde{\D}_R ( e^{sX} ) ( \phi) =   (\phi'(0), \phi'(0)) + (\phi'' (0) , \phi(0)) = X t_1 (\phi) + Y 2 t_2 (\phi)$$
  So we can identify $\D_R (X)$ with $X \otimes t_1 + Y \otimes 2 t_2 \in \hhcm \otimes \bR [\Dplus]$. As $t_1 =1$ we retrieve (\ref{right_coaction_A_on_H}).
   Further, $\tilde{\D}_R ( e^{sY} )(\phi) = \phi \tr ( e^s ,0) = (e^s \phi'(0),\phi(0))$. So
  $$\D_R (Y) (\phi) = {{\frac{d}{ds}}|}_{s=0} \tilde{\D}_R ( e^{sY} ) ( \phi) = (\phi'(0),\phi(0)) = Y t_1 (\phi)$$
  We identify $\D_R ( Y)$ with $Y \otimes t_1 = Y \otimes 1$. 
  So (working with $\field = \bR$) we recover (\ref{right_coaction_A_on_H}), which as we already showed defines a right coaction (for general $\field$).

\subsection{The bicrossproduct $\UCM$}

 We now manufacture a second bicrossproduct $\UCM := \hcmdual$, which we equip with a nondegenerate dual pairing with $\HCM$.
  The action and coaction used to construct the bicrossproduct can again  be motivated by considering the group   factorisation (\ref{D=B+_bowtie_D0}) of $\Diff^+ ( \bR)$.
  First of all:

\begin{lemma}\label{right_module_algebra} $\hhcmdual$ is a right $\ahcmdual$-module algebra via the action defined by
\begin{equation}
\label{right_action_A*_on_H*}
\a \tl z_n =  n \a \b^{n-1},\quad \b \tl z_n =  -\b^n\quad (n \geq 2)
\end{equation}
\end{lemma}
\begin{prf}  It is enough to check that the action (\ref{right_action_A*_on_H*}) defined on generators is compatible with the algebra relations. For example,
$$(\a \tl z_m ) \tl z_n = m ( \a \b^{m-1} ) \tl z_n = m( n-m+1) \a \b^{m+n-2}$$
Hence $\a \tl ( z_m z_n - z_n z_m ) = (n-m) (m+n-1) \a \b^{m+n-2} = \a \tl  [z_m , z_n ]$. \end{prf}

 \begin{lemma}
 \label{lemma_left_coaction_H*_on_A*} 
 $\ahcmdual$ is a left $\hhcmdual$-comodule coalgebra via the coaction  defined on generators by 
  \begin{equation}
  \label{left_coaction_hhcmdual_on_ahcmdual}
 \D_L ( z_n ) = \sum_{j=2}^n {\tiny \nchoosej} \a^{j-1} \b^{n-j} \otimes z_j
 \end{equation}
  and extended to all of $\ahcmdual$ via 
 $\D_L (hg) = \sum ( h^{\coleft} \triangleleft g_{(1)} )  {g_{(2)}}^{\coleft} \otimes h^{\comid} {g_{(2)}}^{\comid}$.
 \end{lemma}
 \begin{prf}  In the same way as Lemma  \ref{lemma_right_coaction_A_on_H} one can check that these formulae extend to a left coaction. In particular it is straightforward  that $( \D \otimes \id) \D_L (z_n) =  (\id \otimes \D_L ) \D_L (z_n)$, and
 $\D_L ( z_m z_n ) - \D_L (z_n z_m ) = (n-m) \D_L ( z_{m+n-1} )$. 
 \end{prf}

 As in Theorem \ref{compatible} it can be checked that the right action (\ref{right_action_A*_on_H*}) and left coaction (\ref{left_coaction_hhcmdual_on_ahcmdual}) are compatible in the sense of \cite{SM_book}, Theorem 6.2.3. Then:

 \begin{prop}
 \label{dual_of_HCM}
 \label{defn_Adual_bicross_Hdual}
 The bicrossproduct Hopf algebra $\UCM :=\hcmdual$  has generators $z_n$ ($n \geq 2$), $\a^{\pm 1}$,  $\b$  with relations $[\a,\b]=0$, 
 $$[ z_m , z_n ] = (n-m) z_{m+n-1},\quad 
[ z_n ,\a ] = -n \a \b^{n-1}, \quad [ z_n , \b ] = \b^n, \quad \D( \a) = \a \otimes \a$$
\begin{equation}\label{defn_UCM}
\D( \b ) = \a \otimes \b + \b \otimes 1,\quad
  \D ( z_n ) = z_n \otimes 1 + \sum_{j=2}^n {\tiny \nchoosej} \a^{j-1} \b^{n-j} \otimes z_j
  \end{equation}
  with antipode and counit defined accordingly.
  \end{prop}
 
 $\UCM$ is an analogue of  the finite-dimensional bicrossproduct $kM \rlbicross k[G]$ of Section \ref{group_factorisations}.
  By the general theory of bicrossproduct Hopf algebras  \cite{SM_book}:
 
 \begin{thm} There is a nondegenerate dual pairing of $\UCM$, $\HCM$, given by 
 \begin{equation}\label{pairing_hcm_hcmdual}
  < {\bf z}  \xi \;, \; {\bf t}  x> := <{\bf z},{\bf t}> < \xi,x>
  \end{equation}
  for all ${\bf z} = z_{m_1} \ldots z_{m_p} \in \ahcmdual$, $\xi = \a^i \b^j \in \hhcmdual$, ${\bf t} = t_{n_1} \ldots t_{n_A} \in \ahcm$, $x = X^r Y^s \in \hhcm$, where on the right hand side we use the pairings defined in (\ref{pairing_A_Adual},\ref{pairing_X^j_Y^k_with_a^p_b^q}).
 \end{thm}

\subsection{Relation to group factorisation II}
\label{relation_to_gp_fact_2}

As in Section \ref{relation_to_gp_fact}, we explain how the right action (\ref{right_action_A*_on_H*}) and left coaction  (\ref{left_coaction_hhcmdual_on_ahcmdual}) used in the construction of  $\UCM$ can be recovered from the factorisation (\ref{D=B+_bowtie_D0}).
 Again this is background and not part of the proof.

Recall that the Lie algebra ${\bf d}_+$ can be represented  as differential operators $z_n = x^n {\frac{d}{dx}}$.
Taking $\field = \bR$, each $z_n$ gives a flow on $\bR$  by solving the ODE $x'(t) = x(t)^n$. 
For $z_0$, $x(t) = x(0) +t$, for $z_1$, $x(t) = x(0) e^t$, and for $z_{n+1}$, with $n \geq 1$,
$x(t) = x(0) [ 1 - n x(0)^n t ]^{-1/n}$. 
The  flows defined by $f_t (x(0)) = x(t)$ are:
\begin{eqnarray}
&&z_0 : f_t (x) = x+t, \quad z_1 : f_t (x) = x e^t\nonumber\\
\label{z_n_flow}
&&z_{n+1}, \; n \geq 1 : \; f_t (x) = x [ 1 - n x^n t ]^{-1/n} = x[ 1 + x^n t] + O(t^2 )
\end{eqnarray}
Obviously these are not defined for all $t$. 
From (\ref{left_action_D0_on_B+}) $\phi \tr (a,b) = ( a \phi' (b) , \phi(b))$, for all $\phi \in \Dplus$, $(a,b) \in \Bplus$.
 We use this to recover the  right action of $\ahcmdual$ on  $\hhcmdual$. 
The flow (\ref{z_n_flow}) corresponding to 
$z_{n+1}$ ($n \geq 1$) is 
$f_\e (x) =  x [ 1 + x^n \e ] + O( \e^2 )$, 
 hence
${\frac{ d f_\e}{dx}} (x) = 1 + (n+1) x^n \e + O( \e^2 )$.
So for $z_{n+1}$, formally we have 
$$f_\e \tr (a,b) = ( a + (n+1) a b^n \e + O( \e^2 ), b + b^{n+1} \e + O( \e^2 ))$$
For $\xi \in \hhcmdual$, define  
$( \xi \tl z_{n+1})(a,b) := {{\frac{d}{d \e}}|}_{\e =0} \xi ( f_\e \tr (a,b) )$, as formally $f_\e = e^{\e z_{n+1}}$. 
Hence
$$ ( \a \tl z_{n+1})(a,b) = {{\frac{d}{d \e}}|}_{\e =0}  a ( 1 + (n+1) b^n \e + O( \e^2)) =  (n+1)  a b^n$$
So $\a \tl z_{n+1} = (n+1)  \a \b^n$ for all $n \geq 1$, and 
the formulae for $\b$ follow in the same way. This explains the motivation for (\ref{right_action_A*_on_H*}), and in Lemma \ref{right_module_algebra} we already proved that this is an action as claimed.

Next, for any group factorisation $X = G \bowtie M$, 
 define a $\field$-linear map 
 $$\tilde{\D}_L : \field M \rto \field [ G, \field M], \quad \tilde{\D}_L (m)(g) = m \tl g$$
 where $\field [ G, \field M]$ is the $\field$-vector space of maps $G \rto \field M$. We now take $\field = \bR$.
 For $z_{n+1}$, the flow (\ref{z_n_flow}) is $f_\e (x) =  x[ 1 + x^n \e ] + O(\e^2 )$.
 Hence
 \begin{eqnarray}
 &(\tilde{\D}_L ( e^{\e z_{n+1} } )(a,b) )(x)&=( f_\e \tr (a,b) )(x) = {\frac{ f_\e (ax+b) - f_\e (b)}{ a {f_\e}' (b)}}\nonumber\\
 &&= {\frac{(ax+b)[ 1 + (ax+b)^n \e ] - b [ 1 + b^n \e ]}{ a [1 + (n+1) b^n \e ]}} + O( \e^2)\nonumber\\
 &&= x + \e \sum_{k=2}^{n+1} {\tiny \nplusonechoosek} a^{k-1} b^{n+1-k} x^k + O( \e^2 )\nonumber
 \end{eqnarray}
   Differentiating with respect to $\e$ and evaluating at $\e =0$ gives a map $\D_L ( z_{n+1}) : \Bplus \rto \ahcmdual$ which we can identify with  (\ref{left_coaction_hhcmdual_on_ahcmdual}).

\subsection{Schr\"odinger action}

Starting with a bicrossproduct Hopf algebra $\H \rlbicross \A$, suppose we have a Hopf algebra $\Adual$ equipped with a nondegenerate dual pairing with $\A$. Then:

\begin{lemma} \cite{SM_book} $\Adual$ is a left $\H \rlbicross \A$-module algebra, via the Schr\"odinger action 
\begin{equation}
\label{left_schrodinger}
(h \otimes a) \tr \phi = \sum h \tr \phi_{(1)} < \phi_{(2)} , a>
\end{equation}
where the left action of $\H$ on $\Adual$ is defined by
$(h \tr \phi)(a) = \phi( a \tl h)$.
\end{lemma}

\begin{cor}\label{schrep} $\hhcm$ is a left $\UCM$-module algebra, via the  Schr\"odinger  action 
$$z_n \tr X =  2Y \d_{n,2}, \; \a \tr X = X,\; \b \tr X =1,\; z_n \tr Y = 0, \; \a \tr Y = Y+1,\; \b \tr Y =0$$
\end{cor}
\begin{prf} 
We have $(z \otimes \xi) \tr x = \sum z \tr x_{(1)} < x_{(2)} , \xi>$, for all $z \in \ahcmdual$, $\xi \in \hhcmdual$, $x \in \hhcm$.
So $(z \otimes 1) \tr X = z \tr X$, defined via
$< z \tr X , h> = < X, h \tl z>$. By (\ref{pairing_X^j_Y^k_with_a^p_b^q}),
$< z_n \tr X , \a^t  \b^r > =<X, (\a^t \tl z_n )  \b^r +  r \a^t  \b^{r-1} (\b \tl z_n ) >$
$=(nt-r) < X, \a^t  \b^{n+r-1}>=  \d_{1,n+r-1} (nt -r) \d_{s,0} = 2 \d_{n,2} \d_{r,0} t$ 
(since $n \geq 2$) whereas $< Y, \a^t   \b^r > = t \d_{r,0}$. 
So $z_n \tr X = 2 Y \d_{n,2}$ as claimed.  The other results follow in the same way.
\end{prf}

There is also a  corresponding dual Schr\"odinger  coaction. First of all, we note that bicrossproducts behave well with respect to dual pairings:

\begin{lemma}\cite{SM_book}
Suppose we are given a bicrossproduct 
$\H \rlbicross \A$, together with Hopf algebras $\Adual$, $\Hdual$ equipped with nondegenerate dual pairings with $\A$, $\H$ respectively.
Suppose further that  $\Adual$ is a right $\Hdual$-comodule coalgebra 
via a coaction $\phi \mapsto   \sum \phi^\comid \otimes \phi^\coright$ 
which is  dual to $\tl : \A \otimes \H \rto \A$ in the sense that 
$$< \phi, a \tl h > = < \sum \phi^\comid \otimes \phi^\coright , a \otimes h>,$$
Then $\Hdual$ is a left $\Adual$-module algebra via the left action  defined by 
$$\tr : \Adual \otimes \Hdual \rto \Hdual, \quad 
(\phi \tr z ) (h) := \sum < \phi, h^\coleft > < z , h^\comid >$$
Furthermore,   this left action and right coaction are compatible in the sense of Theorem \ref{bicrossproduct_thm_1}, 
enabling us to form the bicrossproduct $\Hdual \lrbicross \Adual$.
\end{lemma}

Now consider the linear map 
\begin{equation}
\label{right_coaction_Adual}
\label{right_schrodinger_corepn}
\D_R : \Adual \rto \Adual \otimes \Hdual \otimes \Adual, \quad \phi \mapsto \sum  {\phi_\one}^\comid \otimes {\phi_\one}^\coright \otimes \phi_\two
\end{equation}
It is straightforward to check that  $\Adual$ is a right $\Hdual \lrbicross \Adual$-comodule coalgebra via $\D_R$.
 The Schr\"odinger action  (\ref{left_schrodinger}) and coaction (\ref{right_coaction_Adual}) are dual in the sense: 
$$ \tr : \H \lrbicross \A \otimes \Adual \rto \Adual =  (<.,.> \otimes \id_{\Adual}) ( \id_{\H \lrbicross \A} \otimes \t \circ \D_R)$$
where $\t : \Adual \otimes \Hdual \lrbicross \Adual \rto \Hdual \lrbicross \Adual \otimes \Adual$ is the flip map.  We have:

\begin{lemma}
\label{right_schr_corepn_hcm_on_hhcm}
 The Schr\"odinger  coaction (\ref{right_coaction_Adual}) of $\HCM$ on $\hhcm$ is
$$\D_R (X) = X \otimes 1 + 1 \otimes X + Y \otimes 2 t_2, \quad \D_R (Y) = Y \otimes 1 + 1 \otimes Y$$
\end{lemma}

\section{Scalings and deformations}
\label{section_contractions}

In this section we introduce a natural scale parameter $\l \in \field$ into the  Hopf algebras $\HCM$, $\UCM$ of the previous Section. 
 We first define a family of bicrossproducts $\{ \HCMlambda \}_{\l \in \field}$, with $\HCMlambda \cong \HCM$ for each $\l \neq 0$,  while for $\l =0$ (the so-called classical limit) we obtain a commutative Hopf algebra which can be realised as  functions on the  semidirect product $\bR^2\lcross \Dplus$. 
 We construct a natural quotient Hopf algebra $\clh$ of $\HCMlambda$, which for $\l =0$ corresponds  to  the coordinate algebra of the Heisenberg group. 
 We define a second family $\{ \UCMlambda \}_{\l \in \field}$, again all isomorphic to $\UCM$ for $\l \neq 0$, and find a Hopf subalgebra $\ulh$ with a nondegenerate dual pairing with $\clh$. 
  Finally, by passing to an extended bicrossproduct $\ahcmdual \rlbicross  \FBplus_\l$ we identify the expected classical limits of $\UCMlambda$ and $\ulh$.

\subsection{The deformed Heisenberg bicrossproducts $\HCMlambda$ and $\clh$}
\label{hbg_bicross}

\begin{defn} For each $\l \in k$, we define   $\HCMlambda$ to be  the Hopf algebra with generators $X$, $Y$, $\{ \; t_n \; : \; n = 1,2, \ldots \}$, with $t_1 =1$ and relations
\begin{equation}\label{scaled_relations}
[ Y , X] = \l X, \; [ Y , t_n ] = \l (n-1) t_n,\; [X , t_n ] = \l ( (n+1) {t_{n+1}} - 2 {t_2} {t_n} )
\end{equation}
with coproduct and antipode defined by (\ref{Delta_t_n_first},\ref{antipode_t_n},\ref{new_presentation_HCM}).
\end{defn}

For $\l \neq 0$, the map $\HCM \rto \HCMlambda$ given by 
\begin{equation}\label{HCM_to_HCMlambda}
X \mapsto \l^{-2} X, \quad Y \mapsto \l^{-1} Y, \quad t_n \mapsto \l^{1-n} t_n 
\end{equation}
is a Hopf algebra isomorphism.
 For $\l =0$,   (\ref{scaled_relations}) reduces to a commutative Hopf algebra, denoted $\field [\bR^2\lcross \Dplus ]$, with the generators $X$, $Y$, $t_n$ realisable as functions on the semidirect product  $\bR^2\lcross \Dplus$.

\begin{lemma}\label{hopf_ideal}   Let  $\I$ be the two-sided ideal of $\HCMlambda$ generated by $\{ \; \tt_n  := t_n - t_2^{n-1} \; \}_{n \geq 3}$.\\  Then $\D( \I) \subseteq \HCMlambda \otimes \I + \I \otimes \HCMlambda$ and $\e( \I) = 0$.
\end{lemma} 
\begin{prf} First, $\e( \tt_n) =0 \; \forall \; n$, so $\e( \I) =0$. We have
\begin{eqnarray}
& \D( t_n ) & = \sum_{k=1}^n \sum_{i_1 + \ldots + i_k =n} t_{i_1} \ldots t_{i_k} \otimes t_k\nonumber\\
&&= \sum_{k=1}^n \sum_{i_1 + \ldots + i_k =n} (\tt_{i_1} + t_2^{i_1 -1}) \ldots (\tt_{i_k} + t_2^{i_k -1})  \otimes ( \tt_k + t_2^{k-1}) \nonumber\\
&& = \sum_{k=1}^n \sum_{i_1 + \ldots + i_k =n} t_2^{i_1 + \ldots + i_k - k} \otimes t_2^{k-1}\quad \modulostuff\nonumber\\
&& = \sum_{k=1}^n \sum_{i_1 + \ldots + i_k =n} t_2^{n - k} \otimes t_2^{k-1} = \sum_{k=1}^n \tiny{\nminusonechoosekminusone} t_2^{n-k} \otimes t_2^{k-1}\nonumber\\
&&= \sum_{k=0}^{n-1} \tiny{\nminusonechoosek} t_2^{n-k-1} \otimes t_2^k = ( t_2 \otimes 1 + 1 \otimes t_2 )^{n-1} = \D( t_2^{n-1} )\nonumber
\end{eqnarray}
So $\D( \tt_n ) \subseteq \HCMlambda \otimes \I + \I \otimes \HCMlambda$. Since the $\tt_n$ generate $\I$, the result follows.
\end{prf}

\begin{cor} The quotient bialgebra $\clh : = \HCMlambda/\I$ is in fact a Hopf algebra, generated by $X$, $Y$, $t = 2 t_2$ satisfying 
\begin{eqnarray}
&& [Y, X ] = \l X, \quad [X, t ] =  \half  \l  t^2,\quad \D(X) = 1 \otimes X + X \otimes 1 + Y \otimes t,\nonumber\\
&&[Y, t ] = \l t, \;  \D(Y) = Y \otimes 1 + 1 \otimes Y, \; \D(t) = t \otimes 1 + 1 \otimes t,\; S(Y) = -Y,\nonumber\\
 \label{presentation_clh}
 &&S(X)  = -X + Yt,  \quad S(t) = -t,\quad \e(X) = 0 = \e(Y) = \e(t)
 \end{eqnarray}
\end{cor}
\begin{prf} By Lemma \ref{hopf_ideal} the bialgebra structure of $\HCMlambda$ descends to the quotient, and it is straightforward to check that there is a unique antipode $S$ (defined as shown) that gives $\clh$ the structure of a Hopf algebra.
\end{prf}

  We denote $\hhcm$ with the scaled relation $[Y,X]=\lambda X$
 by $\Ulbplus$. Obviously $\Ulbplus$  and  $U( {\bf b_+})$ are isomorphic Hopf algebras for $\l \neq 0$. Finally, $\kt$ is the commutative unital algebra of polynomials in $t$.
 
\begin{prop}\label{clh_is_bicross} $\clh$  is a bicrossproduct  $\kt \lrbicross \Ulbplus$, via the action 
\begin{equation}\label{left_action_Ulbplus_on_kt}
\tr : \Ulbplus \otimes \kt \rto \kt, \quad X \tr t = \half  \l t^2, \quad Y \tr t = \l t
\end{equation}
and coaction $\D_R : \Ulbplus \rto \Ulbplus \otimes \kt$ defined on generators by
\begin{equation}
\label{defn_kt_coaction}
\D_R (X) = X \otimes 1 + Y \otimes t,\quad \D_R (Y) = Y \otimes 1
\end{equation}
and extended by 
$\D_R (gh) = \sum {g_{(1)}}^{\comid}  h^{\comid} \otimes {g_{(1)}}^{\coright} ( g_{(2)} \triangleright h^{\coright} )$,
\end{prop} 
\begin{prf} It is easy to check that the given coaction is well-defined, and    $\kt$ is a left $\Ulbplus$-module algebra,  $\Ulbplus$ a right $\kt$-comodule coalgebra. 
As in Theorem \ref{compatible} it can be checked that action and coaction are compatible in the sense of Theorem \ref{bicrossproduct_thm_1}.   
So we can construct the bicrossproduct $\kt \lrbicross \Ulbplus$, whose presentation using (\ref{bicross_product_coproduct})  coincides with (\ref{presentation_clh}). 
\end{prf}

The three-dimensional Heisenberg group $\hbg$  is the matrix group 
$$\hbg = \{ \; (a,b,c) := \begin{pmatrix}{1 & a & b\cr 0&1&c\cr 0&0&1}\end{pmatrix} 
\; : \; a,b,c \in \field \; \}$$
If we write the coproduct (\ref{presentation_clh}) of $\clh$ using the  matrix notation
\[ \Delta \begin{pmatrix}{1 & Y & X\cr 0&1&t\cr 0&0&1}\end{pmatrix}=\begin{pmatrix}{1 & Y & X\cr 0&1&t\cr 0&0&1}\end{pmatrix}\tens \begin{pmatrix}{1 & Y & X\cr 0&1&t\cr 0&0&1}\end{pmatrix}\]
  we see that for $\l =0$, $\clh$ is isomorphic to the commutative Hopf algebra generated by the coordinate functions $Y(a,b,c) = a$, $X(a,b,c) = b$, $t(a,b,c) = c$
on $\hbg$. We therefore consider $\clh$ to be  a deformation of the Heisenberg group coordinate algebra.

\subsection{The deformed Heisenberg bicrossproducts $\UCMlambda$ and $\ulh$}

\begin{defn} For each $\l \in \field$, we define $\UCMlambda$ to be the Hopf algebra with generators $z_n$ ($n \geq 2$), $\a$, $\b$ and relations 
$$[\a,\b]=0, \; [ {z_m} , {z_n} ] = (n-m) z_{m+n-1},\; [ {z_n} , \a ] = - \l^{n-1} \, n \, {\a} \b^{n-1},\; [ {z_n} , \b ] =  \l^{n-1} \b^n$$
with coproduct and antipode given by (\ref{defn_UCM}).
\end{defn}

For $\l \neq 0$, there is an isomorphism of Hopf algebras $\UCM \rto \UCMlambda$   defined by 
\begin{equation}\label{scaling_map}
z_n \mapsto \l^{n-1} z_n, \quad \a \mapsto   \a,\quad \b \mapsto \l^2 \b
\end{equation}

\begin{defn}
We define $\ulh$ to be the  Hopf subalgebra of $\UCMlambda$ generated by $z :={z_2}$, $\a$, $\b$. The presentation of $\ulh$ is then: 
\begin{eqnarray}
&&[z , \a ] = -2 \l \a \b, \quad [ z, \b] = \l \b^2, \quad [\a,\beta] =0\nonumber\\
\label{relations_U_lambda_heis}
&&\D( \a ) = \a \otimes \a,\quad \D(\b) =  \b \otimes 1+\a \otimes \b, \quad \D(z) = z \otimes 1 + \a \otimes z
\end{eqnarray}
\end{defn}

 $\ulh$ corresponds to a Heisenberg version of the Planck scale Hopf algebra \cite{mo}.
As before, the $\ulh$ are all isomorphic for $\l \neq 0$.
 Now let $U( z)$ be the unital commutative Hopf algebra generated by $z$, with  $\D(z) = z \otimes 1 + 1 \otimes z$.
 
 \begin{prop} $\ulh$  is a bicrossproduct $U(z) \rlbicross  \hhcmdual$, via the right action 
 $ \tl : \hhcmdual  \otimes U(z) \rto \hhcmdual$ defined by 
$\a \tl z =2 \l \a \b$, $\b \tl z =  -\l \b^2$,
 and  left coaction $\D_L : U(z) \rto \hhcmdual  \otimes U(z)$, $h \mapsto \sum h^\coleft \otimes h^\comid$ defined by  
 $$\D_L (z) = \a \otimes z, \quad \D_L (hg) = \sum ( h^{\coleft} \triangleleft g_{(1)} )  {g_{(2)}}^{\coleft} \otimes h^{\comid} {g_{(2)}}^{\comid} \quad \forall \; h,g \in U(z)$$
\end{prop}
\begin{prf}  It is easily checked that \cite{SM_book}, Theorem 6.2.3 applies.
\end{prf}

\begin{lemma} For $\l \neq 0$, there is a nondegenerate dual pairing of $\clh$ and $\ulh$ given by
$< t^i X^j Y^k , z^p \a^q \b^r > = j ! \, (\l q)^k \, \d_{i,p} \d_{j,r}$.
\end{lemma}
\begin{prf} This follows from (\ref{pairing_hcm_hcmdual}) together with (\ref{HCM_to_HCMlambda}) and (\ref{scaling_map}).
\end{prf}

Obtaining the ``correct" classical limit $\l=0$ of $\ulh$ (in the sense of duality with $\field[\Heis]$) is more subtle, since we would like to obtain the universal enveloping algebra $U({\bf heis})$ of the Heisenberg Lie algebra. 
From the geometric point of view, consider the $\bR$-valued  functions $A$, $\{ \a_t \}_{t \in \bR}$, $\b$  on
$\Bplus$:
  \begin{equation}\label{coordinates_Bplus}
A(a,b) =  \log a, \quad \a_t (a,b) = a^t, \quad \b (a,b) = b
\end{equation}
 We have $\a_{t_1} \a_{t_2} = \a_{t_1 + t_2}$, $\a_0 =1$ and  formally, $\a_t = e^{tA}$.
 To treat the case $\l=0$ we wish to work with $A$ as a generator of $\ulh$ rather than $\a = \a_1$.
 This can be formulated  rigorously in two different ways. 
 One well-known approach is by working over the ring of  formal power series $\field[[\l]]$. A second approach which we now sketch is as follows.
  Motivated by (\ref{coordinates_Bplus}), for any $\l \in \field$ define $\FBplus_\lambda$ as the commutative Hopf algebra (over $\field$) generated by $\{ \a_t \}_{t \in \field}$, $A$, $\b$ with $\a_0 =1$, $\a_{t_1} \a_{t_2} = \a_{t_1 + t_2}$, and
$$\D(\a_t ) = \a_t \otimes \a_t, \quad \D(A) = A \otimes 1 + 1 \otimes A, \quad \D(\b) = \a_\l \otimes \b + \b \otimes 1$$
 Then there is a right action of $\ahcmdual$ on $\FBplus_\l$, and left coaction of $\FBplus_\l$ on $\ahcmdual$ defined by 
 \begin{eqnarray}
 &&\a_t \tl z_n = \l^{n-2} \, nt \, {\a_t} \b^{n-1},\quad A \tl z_n = n \l^{n-2} \b^{n-1}, \quad \b \tl z_n = -\l^{n-1} \b^n\nonumber\\
 &&\D_L ( z_n ) = \sum_{j=2}^n  \l^{n-j} {\tiny{\nchoosej}} \a_{\l (j-1)} \b^{n-j} \otimes {z_j}\nonumber
 \end{eqnarray}
  This action and coaction are compatible in the sense of \cite{SM_book}, Theorem 6.2.3, hence   for each $\l \in \field$ there is a bicrossproduct $\ahcmdual \rlbicross  \FBplus_\l$ containing $\UCMlambda$ as a Hopf subalgebra. Define an extended version of $\ulh$, denoted $\extendedulh$, to be the Hopf subalgebra generated by $z = z_2$, $\a = \a_\l$, $A$ and $\b$. For $\l \neq 0$ this corresponds to $\ulh$ adjoined the primitive element $A$, with $[z,A] = -2\b$, while for $\l =0$ we have $\a = \a_0 = 1$, so $z$, $A$ and $\b$ are primitive with relations $[z,A] = -2\b$, $[z,\b] = 0=[A,\b]$. So for $\l = 0$ then $\extendedulh$ is isomorphic to $U({\bf heis})$.
  Similarly, the correct classical limit of $\UCMlambda$ is the cocommutative Hopf algebra generated by primitive elements $\{ z_n \}_{n \geq 2}$, $A$, $\b$, with $[ {z_m} , {z_n} ] = (n-m) z_{m+n-1}$, $[ z_n , A] = -2 \b \d_{n,2}$, and $[ z_n , \b ] = 0 = [A,\b]$. These remarks are for clarification purposes. We do not use this approach  in the sequel.

Next, the Schr\"odinger action of Corollary \ref{schrep} is compatible with the scaling:

\begin{lemma}\label{contracted_left_schr_repn}
 For $\l \neq 0$, $\ulbplus$ is a left $\ulh$-module algebra via
$$z \tr X=2 \l Y,\quad z \tr Y=0, \quad \a \tr X = X, \quad \a \tr Y = Y + \l, \quad
\b \tr X=\l,\quad \b \tr Y=0$$
\end{lemma}

When $\lambda=0$ there is the known action of the Heisenberg algebra on $\field(X,Y)$ by $z=2Y{\del\over\del X}$, $A={\del\over\del Y}$, $\beta={\del\over\del X}$,
and  the Schr\"odinger action of Lemma \ref{contracted_left_schr_repn} should be thought of as a deformation of this.

We note also that it is immediate that the right Schr\"odinger coaction of Lemma \ref{right_schr_corepn_hcm_on_hhcm} restricts to a right coaction of $\clh$ on $\ulbplus$.

Finally, a typical feature of  bicrossproduct Hopf algebras associated to group factorisations where neither factor group is compact is that the actions have  
singularities \cite{SM_hvn,mo}. These singularities do not appear at the algebraic bicrossproduct level, which is why we have not encountered them in constructing $\clh$ and $\ulh$, but rather when one tries to pass to C*- or von Neumann completions. 
We note that a pair of locally compact quantum groups corresponding to $\ulh$ and $\clh$ was previously  constructed by Vaes \cite{vaes}, Example 3.4, applying the techniques of \cite{bs} to group factorisations $X \cong G \bowtie M$, where both $G$ and $M$ are locally compact.  Explicitly, the correspondence between the generators of the Hopf algebras $\ulh$ and $\clh$, and the (unbounded) operators $A_i$, $B_i$, $C_i$ ($i=1,2$) generating these von Neumann algebras is:
 \begin{eqnarray}
 & \ulh : & \a \mapsto A^2_1, \quad \b \mapsto  A_1 B_1, \quad z \mapsto \l C_1 A^2_1\nonumber\\
 & \clh : & X \mapsto \l B_2, \quad Y \mapsto  \half  \l A_2, \quad t \mapsto 2 C_2\nonumber
 \end{eqnarray}
It is natural to ask whether this could be extended to give  faithful representations of  $\HCM$ and $\UCM$  as (unbounded) operators on some Hilbert space, affiliated to a locally compact quantum group. 
An obstacle is the fact that $\Diff^+ (\bR)$ is not locally compact, nor does it have any interesting locally compact subgroups \cite{goldin}. This question will be pursued elsewhere.

\subsection{Local factorisation of $SL_2(\bR)$}

It is natural to ask for a geometrical picture in terms of a group factorisation linked  to the  $\clh$ and $\ulh$ 
 bicrossproducts, in the same way as $\HCM$ and $\UCM$ were shown to be  linked to the factorisation of $\Diff^+ (\bR)$. In this section we show that the relevant group is $SL_2 (\bR)$, which locally (but not globally) factorises as $SL_2 (\bR) \approx \Bplus \bowtie\bR$. 

Using the   coaction  (\ref{defn_kt_coaction})  define a linear map $f : U(\bplus) \rto U(\bplus)$ by 
$f(x) = \sum x^\comid \; < x^\coright , z>$ 
where the pairing $<.,.> : \kt \otimes U(z) \rto \field$ is  $< t^m ,z^n > = 2^m \d_{m,n}$.
 Then $f(X) = 2Y$, $f(Y) =0$, and further  $f(YX) = 2Y(Y+1)$, $f(XY) = 2 Y^2$, hence $f(YX-XY ) = 2Y = f(X)$. 
 Identifying $X$, $Y$ with the generators of $\bplus$, and $z$ with the generator of the Lie algebra $\br$ of $\bR$,  we have a well-defined  left action of $\br$ on $\bplus$, given by $z \tr x = f(x)$, satisfying 
 \begin{equation}\label{left_lie_action}
 z \tr X = 2Y, \quad z \tr Y = 0
 \end{equation}
 The adjoint of the left action (\ref{left_action_Ulbplus_on_kt}) defines a right action of $U( \bplus)$ on $U(z)$:
 \begin{eqnarray}
&&<t^n , z \tl Y> := < Y \tr t^n , z> = <n t^n, z> = 2n \; \d_{n,1},\nonumber\\
&&< t^n , z \tl X> := < X \tr t^n ,z> = < {\tiny{\half}} n t^{n+1} ,z> = {\tiny{\half}}  n \d_{n+1,1} =0\nonumber
\end{eqnarray}
This gives a right action of $\bplus$ on $\br$, satisfying
\begin{equation}\label{right_lie_action}
z \tl X =0, \quad z \tl Y = z
\end{equation}
It is straightforward to check that $\bplus$ and $\br$ equipped with the actions (\ref{left_lie_action}, \ref{right_lie_action}) are a matched pair of Lie algebras in the sense of  \cite{SM_book}, Definition 8.3.1.
 From this point we take $\field = \bR$. Since we have a matched pair, the $\bR$-vector space ${\bf g} := \bplus \oplus \br$ can be given the structure of a Lie algebra, with Lie bracket
 $$[ x , y ]_{\bf g} : = [x,y]_\bplus, \quad  [z,x]_{\bf g}:=z\ra x+ z\la x, \quad \forall \; x,y \in \bplus$$
 Then 
 $[z,X]_{\bf g}=z\ra X+ z\la X=2Y$, $[z,Y]_{\bf g}=z\ra Y+z\la Y=z$, 
 so $\bg \cong \sltwoR$.
  Since both $\Bplus$ and $\bR$ are simply connected Lie groups, we conclude that
$\Bplus \bowtie\bR \approx \SLtwoR$ 
in so far as the actions on the left hand side exponentiate.
We embed $B_+$ and $\bR$ as  subgroups of $SL_2 (\bR)$ by
\begin{equation}
\label{embed_in_SL2R}
(a,b)=\begin{pmatrix}{a^{1\over 2} & a^{-{1\over 2}}b \cr 0& a^{-{1\over 2}}}\end{pmatrix},\quad (c)=\begin{pmatrix}{1 & 0 \cr -c & 1}\end{pmatrix}
\end{equation}

\begin{lemma} There is a local factorisation $\SLtwoR \approx \Bplus \bowtie\bR$ given by
$$\begin{pmatrix}{a & b \cr c & d}\end{pmatrix} \mapsto (d^{-2},d^{-1}b) (-d^{-1}c), \quad d \neq 0$$
\end{lemma}
\begin{prf} For $d \neq 0$,  
$ \begin{pmatrix}{a & b \cr c & d}\end{pmatrix}=
\begin{pmatrix}{d^{-1}  & b \cr 0 & d }\end{pmatrix}
 \begin{pmatrix}{1 &  0 \cr  d^{-1}c &1}\end{pmatrix}$. Then apply   (\ref{embed_in_SL2R}). 
\end{prf}

The resulting  left action of $\bR$ on $\Bplus$ and right action of $\Bplus$ on $\bR$ are given by:
$$c \tr (a,b) = ( {\frac{a}{(1-bc)^2}} , {\frac{b}{1-bc}}) \quad c \tl (a,b) = {\frac{ac}{1-bc}}$$
As in Sections \ref{relation_to_gp_fact} and \ref{relation_to_gp_fact_2} we can use this factorisation to rederive  the actions and coactions  used to construct the bicrossproducts $\clh$, $\ulh$.

\section{Differential calculi over $\ulbplus$ and $\clh$}
\label{section_fodc}

The model obtained above can be seen as a variant of one of a  family of previously-studied bicrossproducts, which  act on noncommutative algebras denoted $U_\l({\bf b_+^n})$
 of varying dimension $n$.
  First of all, in  \cite{SM:thesis}, a bicrossproduct $U({\bf so_3})\rlbicross \bC[ B_+^3]$  associated to  a local factorisation of $SO(3,1)$ was constructed. This bicrossproduct can be regarded  as  corresponding to a  deformation of the  Euclidean group of motions, and acts naturally on the 
  algebra $U_\l({\bf b_+^3})$. 
  Similarly, in \cite{mr}, a bicrossproduct  $U({\bf so_{3,1}})\rlbicross \bC[B_+^{3,1}]$ associated to a factorisation of $SO(3,2)$ was constructed, and interpreted as   corresponding to a  deformation of the Poincar\'e group. This bicrossproduct acts naturally on the 
  algebra   $U_\l ({\bf b_+^{3,1}})$.
 These and other examples have been widely studied in the mathematical physics literature.
   Our new example coming from the Connes-Moscovici algebra has an analogous geometrical picture.

From this point of view it is natural to  extend the theory to include noncommutative differential geometry both on $\ulbplus$ and on $\clh$ as  `coordinate algebras'.  Recall \cite{woro} that 
 a first order differential calculus (FODC) over an algebra $\A$ is an $\A$-bimodule $\Om^1$ with a linear map $d : \A \rto \Om^1$ such that
 \begin{enumerate}
 \item{$d$ obeys the Leibniz rule $d(ab) = (da) b + a db$, for all $a$, $b \in \A$.}
 \item{$\Om^1$ is the linear span of elements $a db$.}
 \end{enumerate}
 If in addition $\A$ is a right $\H$-comodule algebra for some  Hopf algebra $\H$,  then  $\Om^1$ is said to be right covariant if there exists a linear map 
 $\D_R : \Om^1 \rto \Om^1 \otimes \H$, extending the coaction $\D_R$ of $\H$ on $\A$, in the sense that
 $$\D_R ( a \om b ) = \D_R (a) \D_R (\om) \D_R (b), \quad \D_R ( da) = (d \otimes \id) \D_R (a)$$
 for all $a$, $b \in \H$, $\om \in \Om^1$. Left covariant and bicovariant are defined similarly.
 Covariant FODC over finite-dimensional bicrossproducts have been intensively studied, in particular see \cite{emn}, which deals with the case of bicrossproducts arising from finite group factorisations. 

\subsection{Differential calculi over $\ulbplus$}

There is a standard calculus on $\ulbplus$, 
 the so-called Oeckl calculus \cite{oeckl} generated as a  left $\ulbplus$-module $\Om^1$ by $dX$, $dY$, with 
\begin{equation}
\label{oeckl_calc}
[ dX , X] =0 = [ dX, Y], \quad [ dY, X] = \l dX, \quad [ dY, Y ] = \l dY
\end{equation}
The right Schr\"odinger coaction (\ref{right_schrodinger_corepn}) 
gives $\ulbplus$ the structure of a $\clh$-comodule algebra.
This right coaction is dual to the left Schr\"odinger action (Lemma \ref{contracted_left_schr_repn}) in  that 
$a \tr x = \sum x^\oneb < a, x^\two>$ for all $x \in \ulbplus$, $a \in \ulh$. 
This extends to a right coaction of $\clh$ on $\Om^1$, with
\begin{equation}
\label{coaction_clh_on_ulbplus}
\D_R (dX) = dX \otimes 1 + dY \otimes t, \quad \D_R ( dY) = dY \otimes 1
\end{equation} 

\begin{lemma} The Oeckl calculus (\ref{oeckl_calc}) is covariant under the right coaction (\ref{coaction_clh_on_ulbplus}) if and only if $\l =0$.
\end{lemma}
\begin{prf} $\D_R( (dX)X )-\D_R( X dX) = \l ( dX \otimes t + {\frac{1}{2}} dY \otimes t^2)$, but $[ dX,X] =0$.
\end{prf}

This parallels what is found for the higher dimensional bicrossproducts  mentioned above, where it is known that the natural translation-invariant calculus on $U_\l({\bf b_+^n})$ is not covariant under the bicrossproduct symmetry group.

\begin{thm}
\label{unique_2d_fodc}
 There exists a unique 
  FODC $\Om^1$ over $\ulbplus$ such that:
\begin{enumerate}
\item $\Om^1$ has basis (as a left $\ulbplus$-module) $\{ dX,dY \}$.
\item $\Om^1$ is covariant under the  right coaction (\ref{coaction_clh_on_ulbplus}) of $\clh$. 
\end{enumerate}
The explicit presentation of $\Om^1$ is then
\begin{eqnarray}
&&(dX)X = X dX, \quad (dX) Y = (Y - \l/2 ) dX\nonumber\\
&&(dY)X = (\l/2) dX + X dY, \quad (dY)Y = (Y + \l/2 ) dY\nonumber
\end{eqnarray}
\end{thm}
\begin{prf} It follows from our assumptions that
\begin{eqnarray}
&&(dX) X = a_1 dX + a_2 dY, \quad (dX)Y = b_1 dX + b_2 dY\nonumber\\
&&(dY)X = c_1 dX + c_2 dY, \quad (dY)Y = e_1 dX + e_2 dY\nonumber
\end{eqnarray}
for some $a_1, \ldots , e_2 \in \ulbplus$. Then from $d(YX) - d(XY)= \l dX$, we have $c_1 = b_1 - Y + \l$, $c_2 = X + b_2$. Next, 
\begin{eqnarray}
&&\D_R ( (dY)Y) = (dY) Y \otimes 1 + dY \otimes Y= e_1 dX \otimes 1 + ( e_2 \otimes 1 + 1 \otimes Y) ( dY \otimes 1)\nonumber\\
&&\D_R ( e_1 dX + e_2 dY ) = \D (e_1) ( dX \otimes 1) + [ \D( e_1) (1 \otimes t) + \D( e_2 ) ] ( dY \otimes 1)\nonumber
\end{eqnarray}
Hence $\D(e_1) = e_1 \otimes 1$, $\D(e_2) + \D(e_1) (1 \otimes t) = e_2 \otimes 1 + 1 \otimes Y$, so 
$e_1 = \ep_1$,  $e_2 = Y - \ep_1 t + \ep_2$, where $\ep_1$, $\ep_2$ are constants. 
Applying $\D_R$ to the other expressions gives
\begin{eqnarray}
&&(dX) X = (X+\aap_1) dX + \aap_2 dY, \quad (dX)Y = (Y- \l/2) dX + (\aap_1 /2) dY\nonumber\\
&&(dY)X = (\l/2) dX + (X+ \aap_1/2) dY, \quad (dY)Y = (Y+ \l/2) dY\nonumber
\end{eqnarray}
where $\aap_1$, $\aap_2$ are constants. Then  the constraint 
$(dX)YX - (dX)XY = \l (dX)X$, and similarly for $dY$   gives the result.
\end{prf}

We note further that the left coaction $\ulbplus \rto \clh \otimes \ulbplus$ given by $X \mapsto 1 \otimes X$, $Y \mapsto 1 \otimes Y + Y \otimes 1$ extends to a left coaction on this $\Om^1$ making it into a left covariant FODC. The left and right coactions are compatible, hence $\Om^1$ is in fact bicovariant under $\clh$.

\subsection{Differential calculi over $\clh$}

Similarly, one would like a calculus on the bicrossproduct quantum group itself. In previously-studied    higher dimensional cases it has been found  that any bicovariant calculus needs to have extra non-classical generators. In our case we would expect a four dimensional calculus on $\clh$.
We are therefore  interested to find  FODC over $\clh$ which are bicovariant with respect to the coactions induced from the coproduct. This implies
\begin{equation}\label{left_right_coaction_clh}
\D_L ( dX) = 1 \otimes dX + Y \otimes dt, \quad \D_R ( dX) = dX \otimes 1 + dY \otimes t
\end{equation}
while both $dY$ and $dt$ are left- and right-invariant. 

As shown by Woronowicz \cite{woro}, covariant FODC can be classified in terms of one-sided ideals of the dual Hopf algebra  invariant under the (left or right) adjoint coaction of the dual on itself. 
Using the presentation (\ref{relations_U_lambda_heis}) of $\ulh$ it is straightforward to give a complete list of right-covariant FODC over $\clh$ of dimension at most 4. Then one could check by hand for bicovariance. This would be very laborious and we prefer to proceed directly.
First of all:

\begin{thm}\label{3d_left_cov_fodc}
 Let $\Om^1$ be the $\clh$-bimodule  with left module basis $\{ dX, dY, dt  \}$ and relations
\begin{eqnarray}
&&(dX) X = X dX + \lambda X dt, \quad (dX) Y = Y dX, \quad (dX) t = t dX + \lambda t dt\nonumber\\
&&(dY) X = X dY + \lambda dX, \quad (dY) Y = (Y+\lambda) dY, \quad (dY) t = t dY + \lambda dt\nonumber\\
&&(dt) X = X dt, \quad (dt) Y = Y dt, \quad (dt) t = t dt \nonumber
\end{eqnarray}
Then $\Om^1$ is a left-covariant FODC for the left coaction (\ref{left_right_coaction_clh}). 
\end{thm}

\begin{thm}\label{two_3d_right_cov_fodc}
 For $\l \neq 0$, there exist two nonisomorphic right-covariant three dimensional FODC over $\clh$, with basis $\{ dX, dY, dt \}$, extending the two dimensional FODC of Theorem \ref{unique_2d_fodc}. Explicitly, these are
\begin{eqnarray}
&& (dX) X = X dX, \quad (dX) Y = (Y - \l/2) dX, \quad (dX)t = t dX + gt dt\nonumber\\
&& (dY)X = (\l/2) dY + X dY, \quad (dY)Y = (Y + \l/2) dY, \quad (dY) t = t dY + g dt\nonumber\\
&& (dt)X = (X + (g - \l) t) dt, \quad (dt)Y = (Y + g -\l) dt, \quad (dt) t = t dt\nonumber
\end{eqnarray}
with $g = 0$ or $\l /2$. 
\end{thm}
\begin{prf} We start in the same way as Theorem \ref{unique_2d_fodc}, with the given relations 
$(dX) X = X dX$,...,$(dY)Y = (Y + \l/2) dY$ 
together with 
\begin{eqnarray}
&&(dX)t = c_1 dX + c_2 dY + c_3 dt, \quad (dY)t = g_1 dX + \ldots \nonumber\\
&&(dt) X = h_1 dX + \ldots, \quad (dt) Y = j_1 dX + \ldots, \quad (dt)t = k_1 dX + k_2 dY + k_3 dt\nonumber
\end{eqnarray}
for some $c_1$, ... ,$k_3 \in \clh$. Applying $\D_R$ to both sides of $(dY)Y$, $(dY)t$, $(dt)Y$, $(dt)t$, and using the relation $d(Yt) - d(tY) = \l dt$, we have
$$(dY) Y = \fp_1 dX + (Y - \fp_1 t + \fp_2 ) dY + \fp_3 dt, \quad (dY) t = \gp_1 dX + ( \gp_2 + ( 1 - \gp_1)t ) dY + \gp_3 dt$$
$$(dt) t = \kp_1 dX + (\kp_2 - \kp_1 t) dY + (t+ \kp_3) dt,\quad
 (dt)Y = \gp_1 dX + (\gp_2 - \gp_1 t) dY + (Y + \gp_3 - \l) dt$$
where $\fp_1$, ...,$\kp_3$ are scalars. Doing the same for $(dX)t$ and $(dt)X$ gives 
\begin{eqnarray}
&(dX) t = &[ (1 + \gp_1)t + \l \kp_1 /2 ] dX + [ -\gp_1 t^2 + (\gp_2 - \l \kp_1 /2 )t +(\hp_2 + \l \kp_2 /2) ] dY \nonumber\\
&&+ [ \gp_3 t + (\hp3 + \l \kp_3 /2)] dt\nonumber\\
&(dt)X =&(\gp_1 t) dX + (\l/2)(\kp_2 - \kp_1 t) dY + (\l/2)(2t+ \kp_3) dt\nonumber
\end{eqnarray}
Demanding consistency of all possible relations $(dt)(YX - XY) = \l (dt)X$, ..., $(dY)(Xt - tX) = (\l/2) (dY) t^2$ gives the result. 
\end{prf}

It is straightforward to check that none of the covariant FODC of Theorems \ref{3d_left_cov_fodc} and \ref{two_3d_right_cov_fodc} are bicovariant. In fact:

\begin{thm} Let $\Om^1$ be a three-dimensional FODC over $\clh$, with basis $\{ dX, dY, dt \}$. Then $\Om^1$ cannot be be bicovariant. 
\end{thm}
\begin{prf} This follows in exactly the same way as Theorem \ref{unique_2d_fodc}. 
\end{prf}

We now look for four-dimensional covariant FODC $\Om^1$ over $\clh$.

\begin{thm} Suppose that $\Om^1$ is a four-dimensional right-covariant  FODC $\Om^1$ over $\clh$, with basis (as a left $\clh$-module) $\{ dX, dY, dt, \theta \}$, with $\theta a - a \theta = da$ for all $a \in \clh$, and $\D_R (\theta) = \theta \otimes 1$.
  Such an  $\Om^1$ cannot contain as a sub-bimodule the two dimensional calculus of Theorem \ref{unique_2d_fodc}.
\end{thm}
 \begin{prf}
 As before, write
$$(dX)X = a_1 dX + a_2 dY + a_3 dt + a_4 \theta \quad \ldots \quad (dt)t = k_1 dX + k_2 dY + k_3 dt + k_4 \theta$$
for some $a_1$, ... ,$k_4 \in \clh$.
Applying right-covariance and the relations $[Y,X] = \l X$ and so on gives
\begin{eqnarray}
&(dX) X =& [ X + \fp_1 t^2 + ( \l - 2 \ep_1) t + \aap_1 ] dX + \nonumber\\
&&[ - \fp_1 t^3 + (\fp_2 - \l/2) t^2 + (2 \ep_1 - \aap_1)t + \aap_2 ] dY +\nonumber\\
&&[ \fp_3 t^2 + 2 \ep_3 t + \aap_3 ] dt + [ \fp_4 t^2 + \aap_4 ] \theta\nonumber\\
&(dX)Y =& [ Y + \fp_1 t + (\ep_1 - \l) ] dX +  [ - \fp_1 t^2 + (\fp_2 - \ep_1 ) t + \ep_2 ] dY +\nonumber\\
&& [ \fp_3 t + \ep_3 ] dt + [ \fp_4 t + \ep_4 ] \theta\nonumber\\
&(dX)t =&[(1 + \gp_1)t + ( \cp_1 +(\l/2) \kp_1)] dX +\nonumber\\
&& [ - \gp_1 t^2 + (\gp_2 - \cp_1 - (\l/2) \kp_1)t + (\cp_2 + (\l/2) \kp_2)] dY +\nonumber\\
&& [ \gp_3 t + (\cp_3 +(\l/2) \kp_3)] dt + [ \gp_4 t + (\cp_4 + (\l/2) \kp_4)] \theta\nonumber\\
&(dY)X=& [ \fp_1 t + \ep_1 ] dX + [X - \fp_1 t^2 + (\fp_2 - \ep_1) t + \ep_2 ] dY +\nonumber\\
&& [ \fp_3 t + \ep_3 ] dt + [ \fp_4 t + \ep_4] \theta\nonumber\\
&(dY)Y=& \fp_1 dX + [ Y - \fp_1 t + \fp_2 ] dY + \fp_3 dt + \fp_4 \theta\nonumber\\
&(dY)t=& \gp_1 dX + [\gp_2 + (1- \gp_1)t ] dY + \gp_3 dt + \gp_4 \theta\nonumber\\
&(dt)X=& [ \gp_1 t + \cp_1 ] dX + [ - \gp_1 t^2 + (\gp_2 - \cp_1 ) t + \cp_2 ] dY + \nonumber\\
&& [ X + (\gp_3 - \l) t + \cp_3 ] dt + [ \gp_4 t + \cp_4 ] \theta\nonumber\\
&(dt)Y=& \gp_1 dX + [ \gp_2 - \gp_1 t ] dY + [ Y + \gp_3 - \l ] dt + \gp_4 \theta\nonumber\\
\label{four_dim_right_cov}
&(dt)t=& \kp_1 dX + [ \kp_2 - \kp_1 t] dY + [ t + \kp_3] dt + \kp_4 \theta
\end{eqnarray}
for scalars $\aap_1$, ... , $\kp_4$. There are many  constraints imposed by $(dX)[Y,X] = \l (dX)X$ and so on, we do not list these. 
For $\Om^1$ to  contain as a sub-bimodule the calculus of Theorem \ref{unique_2d_fodc},  we need $(dX)X = X dX$ in (\ref{four_dim_right_cov}), which implies
$$X + \fp_1 t^2 + ( \l - 2 \ep_1) t + \aap_1=X, \quad
 - \fp_1 t^3 + (\fp_2 - \l/2) t^2 + (2 \ep_1 - \aap_1)t + \aap_2=0$$
Hence $\l = 2 \ep_1$, $\aap_1 =0$ and $2 \ep_1 = \aap_1$, which has no solution for $\l \neq 0$.
\end{prf}

\begin{thm} Suppose that $\Om^1$ is a  four-dimensional 
 bicovariant FODC over $\clh$ with basis  $\{ dX, dY, dt, \theta \}$, with 
$\theta a - a \theta = da$ for all $a \in \clh$, 
$\D_L (\theta) = 1 \otimes \theta$ and $\D_R (\theta) = \theta \otimes 1$.
 Then no such $\Om^1$ can exist.
 \end{thm}
\begin{prf} Starting with the relations (\ref{four_dim_right_cov}) and applying $\D_L$ to $(dt)t$, $(dt)Y$, ... ,$(dX)t$ gives $\cp_1 =0$, $\ep_1 = \l$, $\ep_2 = \ep_3 = \ep_4=0$, $\fp_2 = \l$, $\fp_1 = \fp_3 = \fp_4 = 0$, $\gp_3 = \l$, $\gp_1 = \gp_2 = \gp_4 = 0$, $\kp_i =0$ for all $i$.
Hence
$$(dX) X = [X - \l t + \aap_1] dX + [ (\l/2) t^2 + (2 \l - \aap_1)t + \aap_2 ] dY + \aap_3 dt + \aap_4 \theta$$
Applying $\D_L$ to both sides of this gives on the left-hand side:
\begin{eqnarray}
&\D_L ( (dX) X) &= X \otimes dX + 1 \otimes (dX) X + Y \otimes (dt)X + YX \otimes dt\nonumber\\
&&= [ X \otimes 1 + 1 \otimes X - 1 \otimes \l t + \aap_1 (1 \otimes 1) ] (1 \otimes dX) + \other\nonumber
\end{eqnarray}
whereas $\D_L (r.h.s.) = [ \D(X) - \l \D(t) + \aap_1 (1 \otimes 1) ] (1 \otimes dX) + \other$
(where ``other" denotes terms in $dY$, $dt$, $\theta$), and these are inconsistent.
\end{prf}

\section{Appendix: Bicrossproduct description of the left-handed Connes-Moscovici Hopf algebra $\HCMoriginal$}
\label{section_HCM_as_bicross}

We now describe how the original codimension one Connes-Moscovici Hopf algebra $\HCMoriginal$,  defined in \cite{cm}, p206,   is  isomorphic to a right-left bicrossproduct 
  $\hhcm \rlbicross \ahcm$ which we now construct. 
The presentation of $\HCMoriginal$ differs from that of $\HCM$ given in (\ref{HCM_coproduct}) only in that $Y \otimes \d_1$ is replaced by $\d_1 \otimes Y$ in the definition of $\D(X)$. 
The new bicrossproduct is also associated to the factorisation $\Diff^+ (\bR) = \Bplus \bowtie \Dplus$, but in a less natural way than $\HCM$. In particular the  corresponding dual bicrossproduct $\hhcmdual \lrbicross \ahcmdual$ is much more complicated.  This is why we chose to work with $\HCM$ throughout this paper.
Proofs of the assertions in this Section are completely analogous to those given in Section \ref{section_hcm}, and we omit the details. 

\begin{lemma} $\ahcm$ is a right $\hhcm$-module algebra via the action
\begin{equation}\label{right_action_H_on_A}
t_n \tl X = - (n+1) t_{n+1} + 2 t_2 t_n, \quad t_n \tl Y = (1-n) t_n
\end{equation}
 equivalently defined by $\d_n \tl X = - \d_{n+1}$, $\d_n \tl Y = - n \d_n$.
\end{lemma}

\begin{lemma} $\hhcm$ is a left $\ahcm$-comodule coalgebra via the coaction
\begin{equation}\label{left_coaction_A_on_H}
\D_L (X) = 1 \otimes X + 2 t_2 \otimes Y, \quad \D_L (Y) = 1 \otimes Y
\end{equation}
(equivalently,  $\D_L (X) = 1 \otimes X + \d_1 \otimes Y$, $\D_L (Y) = 1 \otimes Y$)
  extended to all of $\hhcm$ via
 $\D_L (hg) = \sum ( h^{\oneb} \triangleleft g_{(1)} )  {g_{(2)}}^{\oneb} \otimes h^{\twob} {g_{(2)}}^{\twob}$.
\end{lemma}

This action and coaction can be derived from (\ref{right_action_B+_on_D0}) as follows. 
 Define a  right action of $\Bplus$ on  $\ahcm$ via
 $$(\xi \tl (a,b))(\phi) := \xi ( \phi \tl (a,b)^{-1} )$$ 
In the same way as  Section \ref{section_hcm} the formulae (\ref{right_action_H_on_A},\ref{left_coaction_A_on_H}) can be recovered. Then:

\begin{prop} The right action (\ref{right_action_H_on_A}) and left coaction (\ref{left_coaction_A_on_H})
 are compatible in the sense of \cite{SM_book}, Theorem 6.2.3.
\end{prop}

This means that there is a well-defined right-left bicrossproduct Hopf algebra $\hhcm \, \rlbicross \, \ahcm$. Using \cite{SM_book}, Theorem 6.2.3 to write out its presentation,  this turns out to coincide  with the presentation of $\HCMoriginal$. Hence:

 \begin{thm}
The   bicrossproduct  $\hhcm \, \rlbicross \, \ahcm$  is isomorphic to the Connes-Moscovici Hopf algebra $\HCMoriginal$.
\end{thm}

\section{Acknowledgements}

We thank Gerald Goldin for very useful discussions about diffeomorphism groups.
We are also grateful for the hospitality of the Perimeter Institute and the Isaac Newton Institute during the time this work was completed. Finally we thank the referee for their helpful comments.


\begin{thebibliography}{99}

\bibitem
{bs} S.~Baaj, G.~Skandalis, {C*-alg\`ebres de Hopf et th\'eorie de Kasparov \'equivariante}, K-theory {\bf 2} (1989) 683-721.

\bibitem
{ck1} A.~Connes, D.~Kreimer, {Renormalization in quantum field theory and the Riemann-Hilbert problem. I. The Hopf algebra structure of graphs and the main theorem}, Comm. Math. Phys. {\bf 210}, no. 1 (2000) 249-273.
 
\bibitem
{cmar} A.~Connes, M.~Marcolli, {$\bQ$-lattices: quantum statistical mechanics and Galois theory}, J. Geom. Phys. {\bf 56}, no. 1 (2006) 2-23.

\bibitem
{cm95} A.~Connes, H.~Moscovici, {The local index formula in noncommutative geometry},  Geom. Funct. Anal. {\bf 5}, no. 2 (1995) 174-243. 

\bibitem
{cm} A.~Connes, H.~Moscovici, {Hopf algebras, cyclic cohomology and the transverse index theorem}, Comm. Math. Phys. {\bf 198} (1998) 199-246.


\bibitem
{cm04a} A.~Connes, H.~Moscovici, {Modular Hecke algebras and their Hopf symmetry}, Mosc. Math. J. {\bf 4}, no. 1, 67-109 (2004) 310.

\bibitem
{emn} J.-P.~Ezin, S.~Majid and F.~Ngakeu,  {Classification of Differentials and Cartan Calculus on Bicrossproducts}, Acta. Applic. Math. {\bf 84} (2004)  193-236.

\bibitem
{fgb} H.~Figueroa, J.~Gracia-Bondia, {Combinatorial Hopf algebras in quantum field theory, I}, 
 Rev. Math. Phys. {\bf 17}, no. 8 (2005) 881-976.


\bibitem
{goldin} G.~Goldin, {Lectures on Diffeomorphism Groups in Quantum Physics}, in Contemporary Problems in Mathematical Physics, Proceedings of the Third International Workshop,
Cotonou, Republic of Benin, 2003, J.~Govaerts, M.~N.~Hounkonnou, A.~Z.~Msezane eds, World Scientific, 2004.
	
\bibitem{kac} G.~I.~Kac, V.~G.~Paljutkin, {Finite ring groups}, Trans. Amer. Math. Soc. {\bf 15} (1966), 251-294.



\bibitem
{lr} J.~Lukierski, H.~Ruegg, {Quantum $\kappa$-Poincar\'e in any dimension},
 Phys. Lett. B {\bf 329}, no. 2-3, (1994) 189-194.

\bibitem
{SM:thesis} S. Majid,  {Non-commutative- geometric Groups by a Bicrossproduct Construction : Hopf Algebras at the Planck Scale}. PhD Thesis, Harvard University Archives, 1988.

\bibitem
{SM_hvn} S.~Majid, {Hopf-von Neumann algebra bicrossproducts, Kac algebra bicrossproducts, and the classical Yang-Baxter equations}, Journal of Functional Analysis {\bf 95} (1991) 291-319.

\bibitem
{SM_qbdg} S.~Majid, {Quantum and braided diffeomorphism groups}, J. Geom. Phys. {\bf 28} (1998) 94-128.

\bibitem
{SM_book} S.~Majid, {Foundations of Quantum Group Theory}, Cambridge University Press, 2000.

\bibitem
{mo} S.~Majid, R.~Oeckl, {Twisting of Quantum Differentials and the Planck Scale Hopf Algebra}, Comm. Math. Phys. {\bf 205} (1999)  617-655.

\bibitem
{mr} S.~Majid, H.~Ruegg, {Bicrossproduct Structure of the $\kappa$-Poincare Group and Non-Commutative Geometry}, Phys. Lett. B. {\bf 334} (1994)  348-354.

\bibitem
{oeckl} R.~Oeckl, {Classification of differential calculi on $U_q(\bplus)$, classical limits, and duality}, J. Math. Phys. {\bf 40}, no. 7 (1999) 3588-3603.


\bibitem{takeuchi} M.~Takeuchi, {Matched pairs of groups and bismash products}, Comm. Algebra {\bf 9} (1981) 841-882.


\bibitem
{vaes} S.~Vaes, {Examples of locally compact quantum groups through the bicrossed product construction}, XIIIth International Congress on Mathematical Physics (London, 2000), 341-348, Int. Press, Boston, MA, 2001.

\bibitem
{woro} S.~L.~Woronowicz,  {Differential calculus on compact matrix pseudogroups (quantum groups)}, Comm. Math. Phys. {\bf 122} (1989) 125-170.




\end{thebibliography}
\end{document}